\documentclass[12 pt]{article}
\usepackage{amssymb,amsmath,amsfonts,amsthm}
\newcommand\BBP{{\mathbb {P}}}
\newcommand\BBI{{\bf 1}}
\newcommand\BBR{{\mathbb {R}}}

\newcommand\BBN{{\mathbb {N}}}
\newcommand\BBZ{{\mathbb {Z}}}
\newcommand\BBE{{\mathbb {E}}}

\newcommand\Var{{\mathrm {Var}}}
\newtheorem {Lemma}{Lemma}[section]
\newtheorem {Theorem}{Theorem}[section]
\newtheorem {Proposition}{Proposition}[section]
\newtheorem {Corollary}{Corollary}[section]

\theoremstyle{definition}
\newtheorem{Definition}{Definition}[section]
\newtheorem{Notation}{Notation}[section]
\newtheorem{Remark}{Remark}[section]
\newtheorem{Example}{Example}[section]

\newcommand\no{\noindent}
\newcommand\ssk{\smallskip}

\newcommand\eps{\varepsilon}
\newcommand\beq{\begin{equation}}
\newcommand\eeq{\end{equation}}

\def\no{\noindent}
\def\ssk{\smallskip}

\def\cov{\mathop{\rm Cov}\limits}

\def\XZ{(X_i)_{i\in{\mathbb Z}}}

\def\EP{(\Omega ,{\cal T},{\mathbb P})}
\def\eps{\varepsilon}\hfuzz=5pt

\begin{document}
\begin{center} {\bf \Large Criteria for Borel-Cantelli lemmas with applications to  Markov chains and dynamical systems}
\vskip15pt

J\'er\^ome Dedecker $^{a}$, Florence Merlev\`{e}de   $^{b}$ {\it
and\/} Emmanuel Rio $^{c}$.
\end{center}
\vskip10pt
$^a$ Universit\'e Paris Descartes, Laboratoire MAP5, UMR 8145 CNRS, 45 rue des Saints-P\`eres, F-75270
Paris cedex 06, France.
E-mail: jerome.dedecker@parisdescartes.fr \\ \\
$^b$ Universit\'e Paris-Est, LAMA (UMR 8050), UPEM, CNRS, UPEC, F-77454 Marne-La-Vall\'ee, France. E-mail: florence.merlevede@u-pem.fr 
\\ \\
$^b$ Universit\'e de Versailles, Laboratoire de math\'ematiques, UMR
8100 CNRS, B\^atiment Fermat, 45 Avenue des Etats-Unis, F-78035 Versailles, France. E-mail: emmanuel.rio@uvsq.fr\vskip10pt


\begin{center}
{\bf Abstract}\vskip10pt
\end{center}

Let $(X_k)$ be a strictly stationary sequence of random variables with values in some Polish space $E$ and common marginal $\mu$, and $(A_k)_{k>0}$ be a sequence of Borel sets in $E$.
In this paper, we give some conditions on   $(X_k)$ and $(A_k)$  under which  the events $\{X_k  \in A_k\}$ satisfy the 
Borel-Cantelli (or strong Borel-Cantelli) property. In particular we prove that, if $\mu (\limsup_n A_n) >0$, the Borel-Cantelli property holds for any  absolutely regular sequence. In case where the $A_k$'s are nested, we show, on some examples, that a rate of  convergence of the mixing coefficients is needed. Finally we give extensions of these results to weaker notions of dependence, yielding applications to 
non-irreducible Markov  chains and dynamical systems.

\section{Introduction}
\label{intro}
\setcounter{equation}{0} 

Let $(\Omega, {\mathcal T}, \BBP)$ be a probability space. Let $\XZ$ be a sequence of random variables defined on $(\Omega, {\mathcal T}, \BBP)$ and with values in some Polish space $E$,  and $(A_k)_{k>0}$ 
be a sequence of Borel sets in $E$. Assume that 
\beq \label{I1}
\BBP ( B_1  )>0  \ \hbox{ and } \sum_{k>0}  \BBP (  B_k ) = \infty, \ \text{where} \ B_k = \{X_k \in A_k\}  . 
\eeq
Our aim in this paper is to find nice sufficient conditions implying the  so-called Borel-Cantelli property
\beq \label{I2}
\sum_{k>0}  \BBI_{B_k}   = \infty   \ \ \hbox{almost surely (a.s.)}   
\eeq
or the stronger one
\beq \label{I3}
\lim_{n \rightarrow \infty}  (S_n /  E_n) = 1 \ \ \hbox{a.s., where }  S_n =  \sum_{k=1}^n  \BBI_{B_k}  \ \text{ and }\  E_n = \BBE ( S_n ) \,  ,
\eeq
 usually called  strong Borel-Cantelli property. The focus will be mainly on irreducible or non-irreducible Markov chains.
Nevertheless we will apply some of our general criteria  to dynamical systems and compare them with the results 
of Kim (2007) and Gou\"ezel (2007) concerning  the transformation defined by  Liverani-Saussol-Vaienti (1999).  
\par\ssk
Let us now recall some known results on this subject. On one hand, if the sequence $\XZ$ is strictly stationary, ergodic, and if $A_k = A_1$ for any positive $k$, then  
$\lim_n n^{  -1} S_n = \mu (A_1)$ a.s., where $\mu$ denotes the law of $X_1$. Hence  (\ref{I2}) holds. 
However, as pointed out for instance by Chernov and Kleinbock (2001), the ergodic theorem cannot be used to handle sequences of sets $(A_k)_k$ such that $\lim_k \mu (A_k) = 0$. 
On the other hand, if the random variables $X_k$ are  independent, then (\ref{I2}) holds  for any sequence $(A_k)_{k>0}$ 
of Borel sets in $E$ satisfying (\ref{I1}) (see Borel (1909), page 252).  Extending this result to non necessarilly independent random variables has been the object of intensive researches.  Let $ {\cal F}_k = \sigma  ( X_i : i \leq k)$ 
and recall that $B_k = \{ X_k \in A_k\}$. L\'evy (1937, p. 249) proved that,  with probability $1$,
\beq
\label{Levy1937}
\sum_{k>0} \BBI_{B_k} = \infty  \ \text{ if and only if }\  \sum_{k>1}  \BBP ( X_k \in A_k \mid {\cal F}_{k-1} ) = \infty.
\eeq
However the second assertion is still  difficult to check in the case of sequences of dependent random variables. 
As far as we know, the first tractable criterion for \eqref{I2} to hold is due to Erd\H{o}s and R\'enyi (1959) and reads as follows: 
\beq \label{CritL2BC}
 \lim_{n\rightarrow \infty}  E_n = \infty \quad \text{and} \quad \lim_{n\rightarrow \infty}  E_n^{-2} \Var (S_n) = 0 \, . 
\eeq
Suppose now that the sequence $B_k = \{ X_k \in A_k \}$ satisfies the following uniform mixing condition: 
\beq \label{I4} 
 | \BBP ( B_k \cap B_{k+n} ) - \BBP (B_k) \BBP ( B_{k+n}) |   \leq \varphi_n \bigl ( \BBP (B_k)  + \BBP ( B_{k+n})  \bigr ) \, .
\eeq
Then, if 
\beq \label{I5}
 \lim_{n\rightarrow \infty}  E_n = \infty \quad \text{and} \quad \sum_{n \geq 1} \varphi_n < \infty \, ,
\eeq
the criterion \eqref{CritL2BC} is satisfied and consequently (\ref{I2}) holds. Furthermore, if (\ref{I5}) holds, then the strong Borel-Cantelli
property (\ref{I3}) also holds, according to Theorem 8 and Remark 7 in Chandra and Ghosal (1998). 
This result has applications to dynamical systems. 
For example, Philipp (1967) considered   the  Gauss map  $T(x) =  1/x $ (mod 1)  and the $\beta$-transforms $T(x) = \beta x$ (mod 1) with $\beta >1$,  with  $(X_k)_{k \geq 0} = (T^k)_{k \geq 0}$ viewed as a random sequence on the probability space $([0,1], \mu)$, where  $\mu$ is the unique $T$-invariant probability measure absolutely 
continuous w.r.t. the Lebesgue measure. For such maps and   sequences  $(A_k)$  of intervals satisfying 
\begin{equation}\label{I1bis} 
\sum_{k>0}  \mu (  A_k ) = \infty,
\end{equation}  he proved that \eqref{I5} is satisfied. 
More recently, Chernov and Kleinbock (2001) proved that \eqref{I5} is satisfied when $(X_k)_{k \geq 0}$ are the iterates of  Anosov diffeomorphisms preserving Gibbs measures and $(A_k)$ belongs to a particular class of rectangles (called EQR rectangles).  We also refer to Conze and Raugi
(2003) for non-irreducible Markov chains satisfying  (\ref{I5}). 

However  some dynamical systems do not satisfy \eqref{I5}. We refer to Haydn {\it et al.} (2013) and Luzia (2014) for examples 
of such dynamical systems and Borel-Cantelli type results, including the strong Borel-Cantelli property.   
In particular, estimates as in  \eqref{I5} are not available for non uniformly expanding maps such as the Liverani-Saussol-Vaienti map 
(1999) with parameter $\gamma \in ]0,1[$.  Actually, for such maps, Kim (2007)  proved in his Proposition 4.2 that for any $\gamma \in ]0,1[$, the sequence of intervals $A_k = [0, k^{1/(\gamma -1)}]$ satisfies \eqref{I1bis} but $(B_k)$ does not satisfy \eqref{I2}. Moreover, there are many irreducible, positively recurrent and aperiodic Markov chains which do not satisfy \eqref{I4} with $\varphi_n \rightarrow 0$ even for regular sets $A_k$, such as the Markov chain considered in Remark \ref{Rioexample} in the case where $A_k = [0,1/k]$ (see Chapter 9 in Rio (2017) for more about irreducible Markov chains). However,  these Markov chains are $\beta$-mixing in the sense of  Volkonski\v{\i} and Rozanov   (1959), and therefore strongly mixing in the sense of Rosenblatt (1956). 

The case where the sequence of events  $(B_k)_{k >0}$  satisfies a strong mixing condition has been considered first by Tasche (1997). 
For $n>0$, let 
\beq
\label{defalphaTasche}
\bar \alpha_n =  \frac 12
\sup \big\{  \BBE  \bigl( \,  |\BBP ( B_{k+n}  \mid {\cal F}_k ) - \BBP ( B_{k+n}  ) | \, \bigr) : k>0 \big\}.
\eeq
Tasche (1997) obtained sufficient conditions for (\ref{I2}) to hold. However these conditions are more restrictive than (\ref{I1}):
even in the case where the sequence $(\bar \alpha_n)_n$ decreases at a geometric rate and $({\mathbb P} (B_k) )_k$ is non-increasing, 
Theorem 2.2 in Tasche (1997) requires the stronger condition $\sum_{k >1}  \BBP(B_k)/\log(k)= \infty$.  Under slower rates of mixing, as a consequence of our Theorem \ref{coralphaBC} (see Remark \ref{remcoralphaBC}), we obtain  that if $( {\mathbb P}(B_k)   )_k$ 
is non-increasing and $\bar \alpha_n \leq C n^{-a} $ for some $a>0$, $ (B_k)_k$ satisfies the  Borel-Cantelli property \eqref{I2} provided that 
\[ 
\sum_{n \geq 1}( {\mathbb P}(B_n) )^{(a+1)/a}= \infty \,  \text{ and } \, \lim_{n \rightarrow + \infty ÷} n^a {\mathbb P}(B_n)  = \infty \, ,
\]
which  improves Item (i) of Theorem 2.2 in Tasche (1997).  Furthermore, we will prove that this result cannot be improved in the specific case of irreducible, positive recurrent and aperiodic Markov chains for some particular sequence  $(A_k)_{k >0}$ of nested sets (see 
Remark \ref{optimalityitem1} and Section \ref{HarrisMC}).  Consequently, for this class of Markov chains, the size property \eqref{I1} is not enough for $(B_k)_{k >0} $ to satisfy \eqref{I2}. 

In the stationary case, denoting by $\mu$ the common marginal distribution, a natural question is then: for sequences of sets $(A_k)_{k >0}$ satisfying the size property \eqref{I1bis}, what conditions could be added  to get the Borel-Cantelli property? Our main result in this direction is Theorem \ref{BCbeta} {\rm (i)} stating that if 
\begin{equation} \label{condbetaBC}
\mu (\limsup_{n } A_n ) >0 \, \text{ and } \,  \lim_{n \rightarrow \infty} \beta_{\infty , 1 } (n) = 0 \, , 
\end{equation}
then $ (B_k)_{k>0}$ satisfies the  Borel-Cantelli property \eqref{I2} without additional conditions on the sizes of the sets $A_k$ 
 (see  \eqref{MR3} for the definition of the coefficients $ \beta_{\infty , 1 } (n) $). 
Notice that the first part of \eqref{condbetaBC} implies  the size property \eqref{I1bis} : this follows from the direct part of the Borel-Cantelli lemma. 
For the  weaker coefficients  $ \tilde \beta_{1 , 1 } (n) $ defined in \eqref{defbetatilde} (resp. $ \tilde \beta_{1 , 1 }^{\rm rev} (n) $ defined in Remark \ref{betarevdef}) and when the $A_k$'s are intervals,  Item (i) of our Theorem \ref{BCbetafaible} implies the Borel-Cantelli property
under the  conditions
\beq
\label{condweakbetaBC}
\mu (\limsup_{n} A_n ) >0 \, \text{ and } \,  \sum_{n >0}   \tilde \beta_{1 , 1 } (n) < \infty \   \Big (\text{resp.}\  \sum_{n >0}  \tilde \beta^{{\rm rev}}_{1,1} (n) < \infty \Big )  .  
\eeq
The proof of this result is based on the following characterization  of sequences $(A_k)$ of intervals satisfying the above condition:
For a sequence $(A_k)$ of intervals, $\mu (\limsup_{n } A_n ) >0$ if and only if there exists a sequence of intervals
$(J_k)$ such that $J_k \subset A_k$ for any positive $k$, $\sum_{k>0} \mu ( J_k) = \infty$ and $(J_k)$ fulfills the asymptotic equirepartition 
property  
\beq
\label{AsymptoticEquirep}
\limsup_n      \left \Vert \frac{ \sum_{k=1}^n   \BBI_{J_k} } {\sum_{k=1}^n \mu  (  J_k)}  \right \Vert_{\infty, \mu} < \infty \, , 
\eeq
where $\Vert \cdot \Vert_{\infty, \mu}$ denotes the supremum norm with respect to $\mu$. 
Up to our knowledge, this elementary result is new. We then prove that, under the mixing condition given in (\ref{condweakbetaBC}), 
the sequence $(\{ X_k \in J_k\})$ has the strong Borel-Cantelli property  (see Item (ii) of Theorem \ref{BCbetafaible}). 
In the case of the  Liverani-Saussol-Vaienti map (1999) with parameter $\gamma \in ]0,1[$, the mixing condition in  (\ref{condweakbetaBC}) holds for $ \tilde \beta^{{\rm rev}}_{1,1} (n) $ and any $\gamma$ in $]0,1/2[$.  
For $\gamma$ in $]0,1/2[$, our result can be applied to prove that $ (B_k)_{k>0}$ satisfies the  Borel-Cantelli property \eqref{I2}  for any sequence $(A_k)$ of intervals satisfying $\mu (\limsup_n A_n ) >0$,  and 
the strong  Borel-Cantelli property \eqref{I3} under the additional condition (\ref{AsymptoticEquirep})  with $J_k= A_k$. 
However, for the LSV map,  Gou\"ezel (2007) obtains the Borel-Cantelli property \eqref{I2} under the  condition
\beq
\label{GouezelCond}
0< \gamma < 1\ \text{ and } \sum_{k>0}  \lambda ( A_k) =  \infty  
\eeq
(but not the strong Borel-Cantelli property).  Now 
$$
\mu ( \limsup_n A_n) > 0 \Rightarrow \lambda ( \limsup_n A_n ) > 0 \Rightarrow \sum_{k>0} \lambda ( A_k) = \infty \, , 
$$
by the direct part of the Borel-Cantelli lemma. Hence, for the LSV map, (\ref{GouezelCond}) is weaker than (\ref{condweakbetaBC}). 
Actually the  condition (\ref{GouezelCond}) is the minimal one to get the Borel-Cantelli property in the case $A_n = [0,a_n]$ (see Example \ref{LSV1} of Section \ref{Examples}). 

A question is then to know if a similar condition to (\ref{GouezelCond})   can be obtained in the setting of irreducible Markov chains.  In this direction, 
we prove that, for aperiodic, irreducible and positively recurrent Markov chains, the renewal measure plays the same role as the Lebesgue measure
for the LSV map. More precisely, if $(X_k)_{k\in \BBN}$ and $\nu$ are respectively the stationary Markov chain and the renewal measure defined in 
Section \ref{HarrisMC}, we obtain the Borel-Cantelli property in Theorem \ref{BCHarris} 
(but not the strong Borel-Cantelli property) for sequences of Borel
sets such that 
\beq
\label{RenewalCrit} 
\sum_{k>0}  \nu ( A_k) = \infty \text{ and }  A_{k+1} \subset A_k \ \text{ for any } k>0 \, ,
\eeq
without additional condition on the rate of mixing. Furthermore we prove in Theorem \ref{BCHarrisconverse} that this condition 
cannot be improved in the nested case. 

\medskip

The paper is organized as follows. In Section \ref{CritBC}, we give some general conditions on a sequence of events $(B_k)_{k>0}$ to satisfy the 
Borel-Cantelli property \eqref{I2}, or some stronger properties (such as the strong Borel-Cantelli property \eqref{I3}). The results of this section, including a more general criterion than 
\eqref{CritL2BC} stated in Proposition \ref{PropCriteriaBC}, will be applied all along the paper to obtain new results in the case where $B_k=\{X_k \in A_k\}$,  under 
various mixing conditions on the sequence $(X_k)_{k >0}$. In  Section \ref{sectionbetaalpha}, we state  our main results for $\beta$-mixing and $\alpha$-mixing sequences;
in Section \ref{sectioncoefffaible}, we consider weaker type of mixing for real-valued random variables, and we give three  examples (LSV map,  auto-regressive processes
with heavy tails and discrete innnovations,  symmetric random walk on the circle)  to which our results apply;
in Section \ref{HarrisMC}, we consider the case where $(X_k)_{k>0}$ is an irreducible, positively recurrent and aperiodic Markov chain: we obtain very precise results, which show in particular that some criteria
of Section \ref{sectionbetaalpha} are optimal in some sense. Section 6 is devoted to the proofs, and some complementary results are given in Appendix (including 
 Borel-Cantelli criteria under pairwise correlation conditions).

\section{Criteria for the Borel-Cantelli properties} \label{sectioncriteriageneral}
\label{CritBC}
\setcounter{equation}{0}
In this section, we give some criteria implying Borel-Cantelli type results. Let $\EP$ be a probability space and $(B_k)_{k>0}$
be a sequence of events. 

\begin{Definition}  \label{BC} The sequence $(B_k)_{k>0}$ is said to be  a Borel-Cantelli sequence in $\EP$ if  
$\BBP  ( \limsup_k B_k) = 1$, or equivalently, $\sum_{k>0}  \BBI_{B_k}   = \infty$ almost surely. 
\end{Definition}

From the first  part of the classical Borel-Cantelli lemma, if $(B_k)_{k>0}$ is a Borel-Cantelli sequence, then $\sum_{k>0} \BBP (B_k) = \infty$.
\par\ssk
 We now define stronger properties. The first one is the convergence in $L^1$. 

\begin{Definition}  \label{L1BC} We say that the sequence $(B_k)_{k>0}$ is a $L^1$ Borel-Cantelli sequence in $\EP$ if 
$\sum_{k>0} \BBP (B_k) = \infty$ 
and
$\lim_{n\rightarrow \infty}  \Vert (S_n /E_n)-1 \Vert_1  = 0$, where $S_n = \sum_{k=1}^n \BBI_{B_k}$ and $E_n = \BBE (S_n)$.
\end{Definition}

Notice that, if $(B_k)_{k>0}$ is a $L^1$ Borel-Cantelli sequence, then $S_n$ converges to $\infty$ in probability as $n$ tends to $\infty$.
Since $(S_n)_n$ is a non-decreasing sequence, it implies that $\lim_n S_n = \infty$ almost surely.  Therefrom $(B_k)_{k>0}$ is a Borel-Cantelli sequence. 
\par\ssk
The second one is the so-called  strong Borel-Cantelli property. 

\begin{Definition}  \label{SBC} With the  notations of Definition \ref{L1BC},  the sequence $(B_k)_{k>0}$ is said to be a strongly Borel-Cantelli sequence if $\sum_{k>0} \BBP (B_k) = \infty$ 
and $\lim_{n\rightarrow \infty}  (S_n/E_n) = 1$ almost surely.   
\end{Definition}

Notice that $\BBE ( S_n/E_n) = 1$. Since the random variables $S_n/E_n$ are nonnegative, by Theorem 3.6, page 32 in Billingsley \cite{Bi1999}, if $(B_n)_{n>0}$ is a strongly 
Borel-Cantelli sequence, then $(S_n/E_n)_{n>0}$ is a uniformly integrable sequence and consequenly $(S_n/E_n)_{n>0}$ converges 
in $L^1$ to $1$. Hence any strongly Borel-Cantelli sequence is a $L^1$ Borel-Cantelli sequence. 

\par\ssk
We start with the following characterizations of the Borel-Cantelli property. 

\begin{Proposition} \label{lmacaracterisationBC}
Let $(A_k)_{k>0}$ be a sequence of events in 
$\EP$  and $\delta \in ]0,1]$ be a real number. The two following statements are equivalent: 
\begin{itemize}
\item[1.] $\BBP  ( \limsup_k A_k)  \geq \delta$. 
\item[2.] There exists  a sequence $(\Gamma_k)_{k>0}$  of events  such that $\Gamma_k \subset A_k$, 
$\sum_{k > 0} \BBP  (  \Gamma_k) =  \infty$ and 
\beq \label{critsequenceGamma}
\limsup_n  \left \Vert \frac{ \sum_{k=1}^n   \BBI_{\Gamma_k} } {\sum_{k=1}^n \BBP  (  \Gamma_k)}  \right \Vert_{\infty} \leq  1/\delta \, .
\eeq
\end{itemize} 
\par
Furthermore, if  there exists a triangular sequence of events $(A_{k,n})_{1 \leq k \leq n}$ with $A_{k,n} \subset A_k$, such that 
${\tilde E}_n := \sum_{k=1}^n \BBP (A_{k,n} ) >0$, $\lim_n {\tilde E}_n = \infty$ and  
$ \bigl( {\tilde E}_n^{-1}  \sum_{k=1}^n \BBI_{A_{k,n}}  \bigr)_{n\geq 1}$ is uniformly integrable, then $\BBP  ( \limsup_k A_k)  > 0$.
\end{Proposition}

Before going further on, we give an immediate application of this proposition which shows that a Borel-Cantelli sequence is characterized by the fact that it contains a subsequence which is a $L^1$ Borel-Cantelli sequence. 

\begin{Corollary} \label{linkBCL1BC}
Let $(A_k)_{k>0}$ be a sequence of events in 
$\EP$  and $\delta \in ]0,1]$ be a real number. Then the following statements are equivalent:
\begin{itemize}
\item[1.] $\BBP  ( \limsup_k A_k)  =1$. 
\item[2.] There exists  a $L^1$ Borel-Cantelli sequence $(\Gamma_k)_{k>0}$  of events  such that $\Gamma_k \subset A_k$. 
\end{itemize} 
\end{Corollary}

Now, if the sets $A_k$  are intervals of the real line, then 
one can construct intervals $\Gamma_k$ satisfying the conditions of  Proposition \ref{lmacaracterisationBC}, as shown by the proposition 
below, which will be applied in Section \ref{Weaktypedep} to the LSV map.  

\begin{Proposition} \label{lmacaracterisationBCint} Let $J$ be an interval of the real line and let $\mu$ be a probability measure on
its Borel $\sigma$-field. 
Let $(I_k)_{k>0}$ be a sequence of subintervals of $J$   and $\delta \in ]0,1]$ be a real number. The two following statements are equivalent:
\begin{itemize}
\item[1.] $\mu  ( \limsup_k I_k)  \geq \delta$. 
\item[2.] There exists  a sequence $(\Gamma_k)_{k>0}$  of intervals  such that $\Gamma_k \subset I_k$, 
$\sum_{k > 0} \mu (  \Gamma_k) =  \infty$ and \eqref{critsequenceGamma} holds true. 
\end{itemize} 
\end{Proposition}

Let us now state some new criteria, which differ from the usual criteria based on pairwise correlation conditions. 
Here it will be necessary to introduce a function $f$ with bounded derivatives up to order $2$. 

\begin{Definition} \label{deffunctionf} Let $f$ be the application from $\BBR$ in $\BBR^+$ defined by $f(x) = x^2/2$ for $x$ in $[-1,1]$ and 
$f(x) = |x| - 1/2$ for $x$ in $]-\infty , -1[ \cup ]1,  +\infty[$. 
\end{Definition}

We now give criteria involving the so defined function $f$. 

\begin{Proposition}  \label{PropCriteriaBC} 
Let $f$ be the real-valued function defined in Definition \ref{deffunctionf} and $(B_k)_{k>0}$ be a sequence of events in 
$\EP$ such that $\BBP (B_1) >0$ and $\sum_{k>0} \BBP (B_k) = \infty$. 
\par\ssk\no
{\rm (i)}   Suppose that there exists a triangular sequence $(g_{j,n})_{1\leq j \leq n}$ of non-negative Borel functions such that 
 $g_{j,n} \leq \BBI_{B_j}$ for any $j$ in $[1,n]$, 
and that  this sequence satisfies the criterion below: if $\tilde S_n = \sum_{k=1}^n g_{k,n} $ and 
$\tilde E_n = \BBE ( \tilde S_n )$,   there exists some increasing sequence $(n_k)_k$ of positive integers such that
\beq \label{CritusualBC}
\lim_{k\rightarrow \infty} \tilde E_{n_k} = \infty \ \text{ and } \lim_{n\rightarrow \infty}  \BBE 
 \bigl(  f \bigl( (\tilde S_{n_k} - \tilde E_{n_k})/  \tilde E_{n_k}   \bigr) \, \bigr)= 0 . 
\eeq
Then $(B_k)_{k>0}$ is a Borel-Cantelli sequence. 
\par\ssk\no
{\rm (ii)}   Let $S_n = \sum_{k=1}^n \BBI_{B_k}$ and $E_n = \BBE (S_n)$.  If 
\beq \label{CritL1BC}
 \lim_{n\rightarrow \infty}  \BBE  \bigl(  f \bigl( (S_n - E_n)/E_n   \bigr) \, \bigr) = 0 ,
\eeq
then $(B_k)_{k>0}$ is a $L^1$ Borel-Cantelli sequence. 
\par\ssk\no
{\rm (iii)}    If 
\beq \label{CritstrongBC}
 \sum_{n>0}  \frac{\BBP (B_n)}{E_n}  \sup_{k\in [1,n]}  \BBE  \bigl(  f \bigl( (S_k-E_k) / E_n  \bigr) \bigr) < \infty ,
\eeq
then $(B_k)_{k>0}$ is a strongly Borel-Cantelli sequence. 
\end{Proposition}

\begin{Remark}  Since $f(x) \leq x^2/2$ for any real $x$, (\ref{CritL1BC}) is implied by the usual $L^2$ criterion \eqref{CritL2BC}, which is the sufficient condition given in Erd\H{o}s and R\'enyi (1959) to prove that $(B_k)_{k>0}$ is a Borel-Cantelli sequence. Moreover,  (\ref{CritstrongBC}) is implied by the more elementary criterion 
\beq \label{CritstrongBCvariance}
 \sum_{n>0}  E_n^{-3} \BBP (B_n) \sup_{k\in [1,n]}    \Var (S_k)  < \infty ,
\eeq
which is a refinement of Corollary 1 in Etemadi (1983) (see also Chandra and Ghosal (1998) for a review).  
\end{Remark}

\section{$\beta$-mixing and $\alpha$-mixing sequences}   \label{sectionbetaalpha}

\setcounter{equation}{0}

In order to state our results, we need to recall the definitions
of the  $\alpha$-mixing, $\beta$-mixing and $\varphi$-mixing coefficients between two
$\sigma$-fields of $\EP$. 

\begin{Definition} The  $\alpha$-mixing coefficient $\alpha ( {\cal A} , {\cal B} )$
between two $\sigma$-fields $\cal A$ and $\cal B$ of $ {\mathcal T}$ is defined by 
\[
2 \alpha (\mathcal{A},\mathcal{B})  =   \sup \{ | \BBE ( |\mathbb{P} (B | {\cal A} )  - \mathbb{P} (B)| ) :   B\in 
\mathcal{B} \}  \, .
\]
One also has $\alpha (\mathcal{A},\mathcal{B})  =   \sup \{ \,  | \BBP  ( A \cap B ) - \BBP (A) \BBP (B) | : (A,B) \in \mathcal{A}\times \mathcal{B}  \}$,
which is the usual definition. 
Now, if $X$ and $Y$ are random variables with values in some Polish space and $\cal A$ and $\cal B$ are the $\sigma$-fields generated respectively by $X$ and $Y$, one can define the $\beta$-mixing coefficient $\beta ( {\cal A} , {\cal B} )$ and the $\varphi$-mixing 
coefficient $\varphi ( {\cal A} , {\cal B} )$ between the $\sigma$-fields $\cal A$ and $\cal B$  by 
$$
\beta ({\cal A}, {\cal B}) =    \BBE \Bigl(   \sup_{B \in {\cal B} }   |\mathbb{P} (B | {\cal A} )  - \mathbb{P} (B)| \Bigr)  
\ \text{ and }\ 
\varphi ({\cal A}, {\cal B}) =    \Big\Vert   \sup_{B \in {\cal B} }   |\mathbb{P} (B | {\cal A} )  - \mathbb{P} (B)| \Big\Vert_\infty \, ,  
$$
where $ \mathbb{P} (\cdot  | {\cal A} ) $ is a regular version of the conditional probability given $ {\cal A}$.  In contrast to the other coefficients $\varphi ({\cal A}, {\cal B}) \not= \varphi ({\cal B} , {\cal A} )$ in the general case.  
\end{Definition}

From these definitions  $2 \alpha ({\cal A}, {\cal B}) \leq \beta ({\cal A}, {\cal B}) \leq \varphi ({\cal A}, {\cal B})  \leq 1$.
According to Bradley (2007), Theorem 4.4, Item (a2), one also has
\beq \label{defequivalentalpha}  
4 \alpha (\mathcal{A},\mathcal{B})= \sup \{\Vert  {\BBE} (Y|%
\mathcal{A})\Vert _{1} : \  Y\text{ $\mathcal{B}$-measurable, } \Vert Y \Vert_\infty = 1 \ \text{and}\  {\BBE} (Y)=0 \} \, .
\eeq

Let us now define the 
the $\beta$-mixing an $\alpha$-mixing coefficients of the sequence $\XZ$. Throughout the sequel 
\beq \label{MR1}
{\cal F}_m = \sigma ( X_k : k\leq m ) \ \text{ and }\ {\cal G}_m = \sigma (  X_i : i \geq m) . 
\eeq
Define the $\beta$-mixing coefficients $\beta_{\infty , 1} (n)$ of $\XZ$ by 
\beq \label{MR3}
\beta_{\infty , 1 } (n) = \beta ( {\cal F}_{-n} , \sigma (X_0) ) \ \hbox{ for any } n>0 \, ,
\eeq
and note that the sequence $(\beta_{\infty, 1}(n))_{n \geq 0}$ is non-increasing.
$\XZ$ is said to be absolutely regular or $\beta$-mixing
if $\lim_{n\uparrow \infty}   \beta_{\infty , 1} (n) = 0$.   Similarly, define the $\alpha$-mixing coefficients $\alpha_{\infty , 1} (n)$  by
\beq \label{defalphabetastrong} 
\alpha_{ \infty,1}(n) = \alpha ({\mathcal F}_{-n}, \sigma(X_0))  \, ,
 \eeq
 and note that the sequence $(\alpha_{\infty, 1}(n))_{n \geq 0}$ is non-increasing.
 $\XZ$ is said to be strongly mixing or $\alpha$-mixing
if $\lim_{n\uparrow \infty}   \alpha_{\infty , 1} (n) = 0$.

 \subsection{Mixing criteria for the Borel-Cantelli properties} 

We start with some criteria when the underlying sequence is $\beta$-mixing and $\mu  ( \limsup_n A_n ) >0$ (see Remark \ref{remarkonthmbeta}). 

\begin{Theorem} \label{BCbeta} Let $(X_i)_{i \in \mathbb Z}$ be a strictly stationary sequence of random variables with values in some Polish 
space $E$. Denote by $\mu$ the common marginal law of the random variables $X_i$. 
Assume that $\lim_{n\uparrow \infty}   \beta_{\infty ,1} (n) = 0$.   Let $(A_k)_{k>0}$ be a sequence of Borel sets in $E$ satisfying $\sum_{k >0} \mu (A_k) = + \infty$. Set $B_k = \{X_k \in A_k\}$ for any positive $k$. 
\par\ssk\no
{\rm (i)} If $\mu  ( \limsup_n A_n ) >0$,  then $(B_k)_{k>0}$ is a Borel-Cantelli sequence.
\par\ssk\no
{\rm (ii)}   Set $E_n = \sum_{k=1}^n \mu (A_k)$ and $H_n = E_n^{-1}  \sum_{k=1}^n \BBI_{A_k} $. 
If  $(H_n)_{n> 0}$ is a uniformly integrable sequence in $(E , {\cal B} (E) , \mu)$, then 
 $(B_k)_{k>0}$ is a $L^1$ Borel-Cantelli sequence. 
\par\ssk\no
{\rm (iii)}   Let $Q_{H_n}$ be the cadlag inverse of the tail function 
$t \mapsto \mu ( H_n>t)$.  Set 
\beq \label{defQ*}
Q^* (0) = 0 \ \text{ and } \ Q^*  (u)  = u^{-1}  \sup_{n>0}  \int_0^u   Q_{H_n} (s) ds \ \text{ for any }  u\in ]0,1]. 
\eeq
If  
\beq \label{critstrongBCbeta}  
\sum_{j>0}  j^{-1}  \beta_{\infty, 1} (j)  Q^*   ( \beta_{\infty, 1} (j) )   <   \infty ,
\eeq
  then $(B_k)_{k>0}$ is a strongly Borel-Cantelli sequence in $\EP$. 
\end{Theorem}

\begin{Remark} \label{remarkonthmbeta}
By the second part of Proposition \ref{lmacaracterisationBC} applied with $A_{k,n} = A_k$, 
if  $(H_n)_{n>0}$ is uniformly integrable, then $\mu  ( \limsup_n A_n ) >0$. Hence {\rm (ii)}  does not apply if $\mu  ( \limsup_n A_n ) =0$. 
On another hand, the map $u \mapsto u Q^* (u)$ is non-decreasing. Thus, if $\beta_{\infty, 1} (j)>0$ for any $j$, 
 (\ref{critstrongBCbeta}) implies that $\lim_{u \downarrow 0} u Q^* (u) = 0$. Then, by Proposition \ref{propquantileuniformint},
$(H_n)_{n>0}$ is uniformly integrable and therefrom $\mu  ( \limsup_n A_n ) >0$. Consequently, if  $\mu  ( \limsup_n A_n )  = 0$,  {\rm (iii)}  cannot be applied if  $\beta_{\infty, 1} (j)>0$ for any $j$.
\end{Remark}

\begin{Remark} If   the sequence $(H_n)_{n>0}$ is bounded in   $L^p ( \mu)$ for some $p$ in 
$]1, \infty]$, $Q^* (u) = {\cal O} ( u^{-1/p} )$ as $u$ tends to $0$. Then, by Proposition \ref{propquantileuniformint},  this sequence is uniformly integrable and consequently, by 
 {\rm (ii)},  $(B_k)_{k>0}$ is a $L^1$ Borel-Cantelli sequence as soon as 
$\lim_{n\uparrow \infty}   \beta_{\infty ,1} (n) = 0$. 
If furthermore $\sum_{j>0} j^{-1} \beta_{\infty, 1} ^{1-1/p}  (j) < \infty$, then, by  {\rm (iii)}, $(B_k)_{k>0}$ is a strongly Borel-Cantelli sequence. 
In particular, if $\mu ( A_i \cap A_j) \leq C \mu (A_i) \mu (A_j)$ 
for any $(i,j)$ with $i\not=j$, for some  constant $C$,  $(H_n)_{n>0}$ is bounded in   $L^2 ( \mu)$, 
and consequently  $(B_k)_{k>0}$ is a strongly Borel-Cantelli sequence as soon as $\sum_{j>0} j^{-1} \sqrt{\beta_{\infty , 1} (j) } < \infty$.  
\end{Remark}

\begin{Remark} Let $S_n = \sum_{k=1}^n\BBI_{A_k} (X_k)$ and $E_n = \BBE(S_n)$. Inequality \eqref{Laststepbeta} 
in the proof of the above theorem applied with $\Gamma_{k,n} = A_k$ gives
\[ \limsup_n
 \BBE  \bigl(  f_n (  S_n - E_n ) \bigr)  \leq 
2 \limsup_n \int_E  G_n  \psi_m  d\mu \, , 
\]
for any $m >0$, where $\psi_m$ is defined in (\ref{intpsi}), $G_n = S_n/E_n$ and $f_n(x) = f(x/E_n)$. It follows that 
\[
\limsup_n
 \BBE  \bigl(  f_n (  S_n - E_n ) \bigr)  \leq 
2 \Vert  \psi_m   \Vert_{\infty}\, 
\]
for any positive integer $m$. Now, from inequality (\ref{intpsi}) in the proof of Theorem \ref{BCbeta}, we have $\Vert \psi_m \Vert_\infty \leq \varphi ( \sigma (X_0) , {\cal F}_{-m} )$. 
Hence, if  $\varphi ( \sigma (X_0) , {\cal F}_{-m} )$ converges to $0$ as $m$ tends to $\infty$, 
then $\lim_n \BBE  \bigl(  f_n ( S_n - E_n ) \bigr)  = 0$ and consequently $(B_k)_{k>0}$ is a $L^1$ Borel-Cantelli sequence (see Item (ii) of Proposition \ref{PropCriteriaBC}). Similarly, one can prove that, if $\varphi (  \sigma (X_0),  {\cal G}_m \, )$ converges to $0$ as $m$ tends to $\infty$, then $(B_k)_{k>0}$ is a $L^1$ Borel-Cantelli sequence.  For other results in the $\varphi$-mixing  setting, see Chapter 1 in Iosifescu and Theodorescu (1969). 
\end{Remark}

Let us now turn to the general case where $\mu  ( \limsup_n A_n )$ is not necessarily  positive. In this case, assuming absolute regularity does not yield any improvement compared to the strong mixing case (see Remark \ref{optimalityitem1} after Corollary \ref{coralphadirect}). 
Below, we shall use the following definition of the inverse function associated with some non-increasing sequence of reals.

\begin{Definition} \label{definversevn} For any non-increasing sequence $(v_n)_{n \in \BBN}$ of reals, the function $v^{-1}$ is defined by 
$v^{-1} (u) = \inf \{  n \in \BBN \, : \, v_n \leq u   \} = \sum_{n \geq 0} \BBI_{\{ u < v_n \}}$.
\end{Definition}

\begin{Theorem} \label{coralphaBC}
Let $\XZ$ be a strictly
stationary sequence of random variables with values in some Polish space $E$. Let $(\alpha_{ \infty,1}(n))_{n \geq 0}$ be its associated sequence of  strong-mixing coefficients defined by \eqref{defalphabetastrong}. 
Denote by $\mu$ the law of $X_0$. Let $(A_k)_{k>0}$ 
be a sequence of Borel sets in $E$ satisfying $\sum_{k >0} \mu (A_k) = + \infty$.  Set $B_k = \{X_k \in A_k\}$ for any positive $k$. 
 Assume that there exist $n_0 >0$, $C>0$, $\delta >0$ and a non-increasing sequence  $(\alpha_*(n))_{n \geq 0}$ such that  for all $n \geq n_0 $, 
\beq  \label{condregulariryalpha}
\alpha_{ \infty,1}(n)\leq C \alpha_*(n) \,  \text{ and } \,  \alpha_{*} (2n) \leq (1-\delta ) \alpha_*(n)  \, .
\eeq 
Suppose in addition that $(\mu (A_n))_{n \geq 1}$ is a non-increasing sequence, 
\beq \label{taschealpha}
\frac{\mu (A_n) }{\alpha_*(n)} \rightarrow \infty \, \text{ as $n \rightarrow \infty$, and } \, \sum_{n \geq 1} \frac{\mu (A_n) }{\alpha_*^{-1} (\mu (A_n))} = \infty \, .
\eeq 
Then $(B_k)_{k>0}$ is a Borel-Cantelli sequence.  
\end{Theorem}
\begin{Remark} \label{remcoralphaBC} Let us first notice that  Theorem \ref{coralphaBC}  still holds with ${\bar \alpha}_n$ defined in \eqref{defalphaTasche} instead of $\alpha_{ \infty,1}(n)$ (the proof is unchanged). To compare Theorem \ref{coralphaBC} with Theorem 2.2 (i) in Tasche (1997), let us consider 
\[
\mu (A_n) \sim C_1 n^{-(r+1)/ (r+2) } (\log n)^{-b} \text{ and } {\bar \alpha}_n \sim C_2 n^{-(r+1) } (\log n)^{-a}
\]
with $r \geq  -1$. Theorem 2.2 (i) in Tasche (1997) requires $ a >1$ and $ b \leq 1$  whereas an application of Theorem \ref{coralphaBC}  gives the weaker  conditions: $(r+2) b \leq a+ r+1$ if $r>-1$ and $a>b$ if $r=-1$. 
\end{Remark}

\begin{Theorem} \label{coralphaSBC}
Let $\XZ$ be a strictly
stationary sequence of random variables with values in some Polish space $E$. Let $(\alpha_{ \infty,1}(n))_{n \geq 0}$ be its associated sequence of  strong-mixing coefficients defined by \eqref{defalphabetastrong}. 
Denote by $\mu$ the law of $X_0$. Let $(A_k)_{k>0}$ 
be a sequence of Borel sets in $E$ satisfying $\sum_{k >0} \mu (A_k) = + \infty$.  Set $B_k = \{X_k \in A_k\}$ for any positive $k$.  
Let $E_n = \sum_{k=1}^n \mu(A_k)$. 
\begin{itemize}
\item[1.] Let $\eta (x) = x^{-1} \alpha_{\infty,1} ([x])$. Assume that $\lim_n E_n^{-1}\eta^{-1} (1/n)  = 0$.  Then 
$(B_k)_{k>0}$ is a $L^1$ Borel-Cantelli sequence. 
\item[2.] Assume that there exist a sequence $(u_n)_{n >0}$ of positive reals such that 
\beq \label{condforalphaSBC}
\sum_{n>0}\frac{\mu(A_n)}{E_n} u_n < \infty  \, \text{ and } \,  \sum_{n>0}\frac{\mu(A_n)}{E^2_n} \alpha_{\infty,1}^{-1} ( E_n u_n/n)  < \infty \, .
\eeq 
Then $(B_k)_{k>0}$ is a strongly Borel-Cantelli sequence. 
\end{itemize}
\end{Theorem}

We now apply these results to rates of mixing  ${\cal O} ( n^{-a} )$ for some positive constant $a$. 

\begin{Corollary} \label{coralphadirect}
Let $(A_k)_{k>0}$ 
be a sequence of Borel sets in $E$ satisfying $\sum_{k >0} \mu (A_k) = + \infty$. For any $k>0$, let $B_k = \{X_k \in A_k\}$. 
Assume that there exists $a >0$ such that $\alpha_{\infty,1}(n) \leq C n^{-a}$, for $n \geq 1$. 
\begin{itemize}
\item[1.] If $\sum_{n \geq 1}(  \mu (A_n) )^{(a+1)/a}= \infty$, $\lim_n n^a \mu (A_n) = \infty$ and $(\mu (A_n) )_{n \geq 1}$ is non-increasing, then $ (B_k)_{k>0}$ is a  Borel-Cantelli sequence. 
\item[2.] If $\lim_n n^{-1/(a+1)} E_n = \infty$ then  $(B_k)_{k>0}$ is a $L^1$ Borel-Cantelli sequence. 
\item[3.] If $\sum_{n>0}n^{1/(a+1)}\mu(A_n) E_n^{-2} < \infty$ then  $(B_k)_{k>0}$ is a strongly Borel-Cantelli sequence. 
\end{itemize}
\end{Corollary}

\begin{Remark} \label{optimalityitem1} According to the second item of  Remark \ref{Rioexample}, Item 1. of Corollary \ref{coralphadirect} 
cannot be improved, even in the $\beta$-mixing  case. 
\end{Remark}

\begin{Remark} \label{remreversedtime} Theorems \ref{coralphaBC} and \ref{coralphaSBC} (and therefore Corollary \ref{coralphadirect}) also hold if the coefficients $\alpha_{\infty,1}(n)$ are replaced by the reversed ones $\alpha_{1,\infty}(n) = \alpha ( \sigma (X_0), {\mathcal G}_n )$ (see Section \ref{sectionremreversedtime} for a short proof of this remark). 
\end{Remark}

\begin{Remark} \label{RemarkcomparisonwithpropBC}
Let $\alpha_{1,1} (n) = \alpha ( \sigma(X_0), \sigma(X_n))$.  From the criteria based on pairwise correlation conditions stated in Annex \ref{pairwisecoeff}, if $\alpha_{1,1} (n) = {\cal O} (n^{-a}) $ with $ a>1$  then  $(B_k)_{k>0}$ is a $L^1$ Borel-Cantelli sequence if $\lim_n n^{-1/(a+1)} E_n = \infty$ (see Remark \ref{RemarkpropBC}), which is the same condition as
in Corollary \ref{coralphadirect}. Now if $\alpha_{1,1} (n)  = {\cal O} (n^{-a}) $ with $ a \in ]0,1[$,  $(B_k)_{k>0}$ is a $L^1$ Borel-Cantelli sequence when $\lim_n n^{- 1 + a/2 } E_n = \infty$  (see Remark \ref{RemarkpropBC}), which is  more restrictive.
Recall  that, for Markov chains $\alpha_{\infty , 1} (n) = \alpha_{1,1} (n)$. Hence criteria based on pairwise correlation conditions 
are less efficient in the context of $\alpha$-mixing Markov chains and slow rates of $\alpha$-mixing.  
\end{Remark}

\section{Weakening the type of dependence} \label{sectioncoefffaible}  
\label{Weaktypedep}

\setcounter{equation}{0}

In this section, we consider stationary sequences of real-valued random variables. In order to get more examples than $\alpha$-mixing or $\beta$-mixing sequences, 
we shall use less restrictive coefficients, where the test functions are indicators of half lines instead of indicators of Borel sets. 
Some exemples of slowly mixing dynamical systems and non-irreducible Markov chains  to which our results apply will be given in Subsection \ref{Examples}. 

\subsection{Definition of the coefficients}

\begin{Definition} The  coefficients $ \tilde \alpha ( {\cal A} , X )$ and $\tilde \beta ( {\cal A} , X) $ 
between a $\sigma$-field $\cal A$ and a real-valued random variable $X$ are defined by 
$$
\tilde \alpha ({\cal A}, X) =   \sup_{t \in {\mathbb R}}  \left \|  {\mathbb E}(\BBI_{X\leq t}|{\mathcal A}) -{\mathbb P}(X \leq t) \right \|_1 
\ \text{ and }\ 
\tilde \beta ({\cal A}, X) =     
\Big\| \sup_{t \in {\mathbb R}} \left | {\mathbb E}(\BBI_{X\leq t}|{\mathcal A}) -{\mathbb P}(X \leq t) \right | \Big \|_1 \, .
$$
The coefficient $\tilde \varphi ( {\cal A} , X) $ between 
 $\cal A$ and $X$ is defined by 
 $$
\tilde \varphi ({\cal A}, X) =   \sup_{t \in {\mathbb R}}  \left \|  {\mathbb E}(\BBI_{X\leq t}|{\mathcal A}) -{\mathbb P}(X \leq t) \right \|_\infty \, .
$$
\end{Definition}
From this definition it is clear  that $ \tilde \alpha ({\cal A}, X) \leq \tilde  \beta ({\cal A}, X) \leq  \tilde  \varphi ({\cal A}, X) \leq  1$.
\par
Let $\XZ$ be a stationary sequence of real-valued random variables. We now define the dependence coefficients of 
$\XZ$ used in this section. 
The  coefficients $ \tilde \alpha_{\infty , 1} (n)$  are defined by 
\beq \label{defalphatilde} 
\tilde \alpha_{ \infty,1}(n) =  \tilde \alpha ({\mathcal F}_{0}, X_n)  \ \hbox{ for any } n>0. 
 \eeq
Here ${\cal F}_0 = \sigma ( X_k : k\leq 0 ) $ (see \eqref{MR1}). 
The  coefficients $ \tilde \beta_{1 , 1} (n)$  and $ \tilde \varphi_{1 , 1} (n)$  are defined by
\begin{equation} \label{defbetatilde}
 \tilde \beta_{1 , 1 } (n)  =  \tilde \beta ( \sigma(X_{0}) , X_n ) \quad \text{and} \quad  \tilde \varphi_{1 , 1 } (n)  =  \tilde \varphi ( \sigma(X_{0}) , X_n ) \hbox{ for any } n>0.
\end{equation}

\subsection{Results}

 \begin{Theorem} \label{BCbetafaible} Let $(X_i)_{i \in \mathbb Z}$ be a strictly stationary sequence of real-valued random variables. 
 Denote by $\mu$ the common marginal law of the random variables $X_i$. 
  Let $(I_k)_{k>0}$ be a sequence of intervals such that 
$\mu (I_1) >0$ and $\sum_{k>0} \mu (I_k) = \infty$. 
Set $B_k = \{X_k \in I_k\}$ for any positive $k$, and  $E_n = \sum_{k=1}^n \mu (I_k)$.
\par\ssk\no
{\rm (i)}  If  $\mu  ( \limsup_n I_n ) >0$  and  
$\sum_{k >0}  \tilde \beta_{1,1} (k) < \infty$, then $(B_k)_{k>0}$ is a Borel-Cantelli sequence.
\par\ssk\no
{\rm (ii)}   
Let $p \in [1, \infty)$ and $q$ be the conjugate exponent of $p$.  If 
$$
   \lim_{n \rightarrow \infty}  \frac{1}{E_n^p} \sum_{k=1}^{n -1} k^{p-1} \tilde \beta_{1,1} (k) =0  \quad \text{and} \quad \sup_{n >0} 
\frac{1} {E_n}   \Bigl \| \sum_{k=1}^n \BBI_{I_k} (X_0) \Bigr \|_q < \infty \, ,
$$ 
then  $(B_k)_{k>0}$ is a $L^1$ Borel-Cantelli sequence. 
\par\ssk\no
{\rm (iii)}    
Let $p \in [1, \infty)$ and $q$ be the conjugate exponent of $p$.  If 
$$
   \sum_{n>0}  \frac{\mu(I_n)}{E_n^2}  \left (\sum_{k=1}^{n-1} k^{p-1}   \tilde \beta_{1,1} (k) \right)^{1/p}   < \infty \quad \text{and} 
\quad \sup_{n >0} 
\frac{1} {E_n}   \Bigl \| \sum_{k=1}^n \BBI_{I_k} (X_0) \Bigr \|_q < \infty \, ,
$$  
then $(B_k)_{k>0}$ is a strongly Borel-Cantelli sequence. 
\par\ssk\no
{\rm (iv)}   
If $\, \lim_{n \rightarrow \infty}  E_n^{-1} \sum_{k=1}^{n -1}  \tilde \varphi_{1,1} (k) =0$, 
   then 
 $(B_k)_{k>0}$ is a $L^1$ Borel-Cantelli sequence. 
\par\ssk\no
{\rm (v)}   
 If $\, \sum_{n>0}  E_n^{-1} \tilde \varphi_{1,1} (n)< \infty$, 
 then $(B_k)_{k>0}$ is a strongly Borel-Cantelli sequence.
\end{Theorem}

\begin{Remark} Item {\rm (v)} on the uniform mixing case can be derived from Theorem 8 and Remark 7 in Chandra and Ghosal (1998).
Note that, if $p=1$, the condition in Item {\rm (iii)} becomes 
$$
\sum_{n>0}   \frac{ \tilde \beta_{1,1} (n)}{E_n}    < \infty \quad \text{and} \quad \sup_{n >0} 
\frac{1} {E_n}   \left \| \sum_{k=1}^n \BBI_{I_k} (X_0) \right \|_\infty < \infty \, .
$$
Note that, for intervals  $(I_k)_{k >0}$ satisfying the condition on right hand, we get the same condition as in {\rm (v)}, but for
$\tilde \beta_{1,1} (n)$ instead of $\tilde \varphi_{1,1} (n)$.
\end{Remark} 

\begin{Remark} \label{betarevdef}
Theorem  \ref{BCbetafaible} remains true if we replace the coefficients $\tilde \beta_{1,1} (n)$ (resp. $\tilde \varphi_{1,1} (n)$) by 
$\tilde \beta^{{\rm rev}}_{1,1} (n)= \tilde \beta ( \sigma(X_{n}) , X_0 )$ (resp. $ \varphi^{{\rm rev}}_{1,1} (n)= \tilde \varphi ( \sigma(X_{n}) , X_0 )$).
\end{Remark}

\begin{Remark}  {\sl Comparison with usual pairwise correlation criteria.\/}  Let us compare Theorem \ref{BCbetafaible} wit the results stated in Annex  \ref{pairwisecoeff} in the case $\mu ( \limsup_n I_n) >0$. From the definition of the coefficients $\tilde \beta_{1,1} (n)$, 
\[
| \BBP ( B_k \cap B_{k+n} ) - \BBP ( B_k ) \BBP ( B_{k+n} ) | \leq \tilde \beta_{1,1} (n) . 
\]
Hence the assumptions of Proposition  \ref{propPC}  hold true with $\gamma_n = \varphi_n = 0$ and $\alpha_n = \tilde \beta_{1,1} (n)$. 
In particular, from Proposition \ref{propPC}{\rm (i)},  if 
\beq  \label{CritAnnexB}
\lim_n E_n^{-2} \sum_{k=1}^n  \sum_{j=1}^k \min ( \tilde \beta_{1,1} (j)  , \mu (I_k) )  = 0 , 
\eeq
  $(B_k)_{k>0}$ is a Borel-Cantelli sequence. For example, if 
$\tilde \beta_{1,1} (n) = {\cal O} ( n^{-a} )$ for some constant $a>1$, then, from Remark \ref{RemarkpropBC}, 
(\ref{CritAnnexB}) holds if $\lim_n n^{-1/(a+1)} E_n = \infty$. In contrast  Theorem \ref{BCbetafaible}{\rm (i)} ensures that $(B_k)_{k>0}$ is Borel-Cantelli sequence as soon as $\sum_{k>0} \tilde \beta_{1,1}(k) < \infty$, without  conditions on the sizes of the intervals $I_k$. Next,
if $\tilde \beta_{1,1} (n) = {\cal O} ( n^{-a} )$ for some $a<1$, then, according to Remark \ref{RemarkpropBC}, (\ref{CritAnnexB}) is fulfilled if $\lim_n n^{-1 + (a/2)} E_n = \infty$. Under the same condition, Theorem \ref{BCbetafaible}{\rm (ii)} ensures that $(B_k)_{k>0}$ is a Borel-Cantelli sequence
if, for some real $q$ in $(1, \infty]$,  
\beq
\lim_n n^{-1 + (a/p)} E_n = \infty  \quad \text{and} \quad \sup_{n >0} 
\frac{1} {E_n}   \Bigl \| \sum_{k=1}^n \BBI_{I_k} (X_0) \Bigr \|_q < \infty \, ,
\eeq
where $p= q/(q-1)$. Consequently Theorem \ref{BCbetafaible}{\rm (ii)} provides a weaker condition on the sizes of the intervals $I_k$ if 
the sequence $ ( \sum_{k=1}^n \BBI_{I_k} (X_0) / E_n)_{n>0}$ is bounded in $L^q$ for some $q>2$. 
\end{Remark}

As quoted in Remark  \ref{remarkonthmbeta}, if $\mu  ( \limsup_n I_n ) =0$ then {\rm (i), (ii), (iii)} of Theorem \ref{BCbetafaible} 
cannot be applied. Instead, the analogue of Theorems \ref{coralphaBC} and \ref{coralphaSBC} and of Corollary \ref{coralphadirect} hold (the proofs are unchanged). 

\begin{Theorem} \label{thmalphafaible}
Let $(X_i)_{i \in \mathbb Z}$ be a strictly stationary sequence of real-valued random variables. 
 Denote by $\mu$ the common marginal law of the random variables $X_i$. 
  Let $(I_k)_{k>0}$ be a sequence of intervals such that 
$\mu (I_1) >0$ and $\sum_{k>0} \mu (I_k) = \infty$. 
Set $B_k = \{X_k \in I_k\}$ for any positive $k$, and  $E_n = \sum_{k=1}^n \mu (I_k)$.
Then the conclusion of Theorem \ref{coralphaBC} (resp. Theorem \ref{coralphaSBC}, Corollary \ref{coralphadirect}) holds by replacing the conditions on
$(\alpha_{\infty, 1}(n))_{n>0}$ and $(A_k)_{k>0}$ in Theorem \ref{coralphaBC} (resp. Theorem \ref{coralphaSBC}, Corollary \ref{coralphadirect}) by the same conditions
on $( \tilde \alpha_{\infty, 1}(n))_{n>0}$ and $(I_k)_{k>0}$. 
\end{Theorem}
\begin{Remark} \label{rmkthmalphafaible}
Theorem \ref{thmalphafaible} remains true if we replace the coefficients $\tilde \alpha_{\infty, 1}(n)$ by $\tilde \alpha_{ 1, \infty}(n) = \tilde \alpha ( {\mathcal G}_n, X_0 )$ where ${\mathcal G}_n = \sigma ( X_i, i \geq n ) $ (see the arguments given in the proof of Remark \ref{remreversedtime}). 
\end{Remark}

\subsection{Examples} \label{Examples}

\begin{Example} \label{LSV1}
Let us consider the so-called LSV map  (Liverani, Saussol and Vaienti  (1999)) defined as follows: 
\begin{equation}\label{int}
\mbox{for $0 < \gamma <1$}, \quad  \theta(x)=
  \begin{cases}
  x(1+ 2^\gamma x^\gamma) \quad  \text{ if $x \in [0, 1/2[$}\\
  2x-1 \quad \quad \quad \ \  \text{if $x \in [1/2, 1]$} \, .
  \end{cases}
\end{equation}
Recall that if $\gamma \in ]0, 1[$, there is only one absolutely 
continuous invariant probability $\mu$ whose density $h$  satisfies $0<c\leq h(x)/ x^{-\gamma}\leq C<\infty$.  Moreover, it has been proved in \cite{DDT}, that the $\tilde \beta^{{\rm rev}}_{1, 1} (n)$ coefficients of weak dependence associated with $(\theta^n)_{n \geq 0}$, viewed as a random sequence defined on $([0,1], \mu ) $,  satisfy $  \tilde \beta^{{\rm rev}}_{1, 1} (n) \leq  \kappa n^{- ( 1- \gamma)/\gamma }$ for any $n \geq 1$ and  some  $\kappa >0 $.  

Let us first recall Theorem 1.1 of Gou\"ezel (2007): let $\lambda$ be the Lebesgue measure over $[0,1]$ and let  $(I_k)_{k>0}$  be a sequence of intervals such 
that 
\begin{equation}\label{Gou}
 \sum_{k >0} \lambda (I_k) = \infty \, .
 \end{equation}
  Then  $B_n = \{\theta^n \in I_n\} $ is a Borel-Cantelli sequence.  If furthermore the intervals $I_k$ are included in $[1/2, 1]$ then $B_n = \{\theta^n \in I_n\} $ is a strongly Borel-Cantelli sequence (this follows from inequality (1.3) in \cite{Go}, and Item (ii) of Proposition   
  \ref{propPC}.) If  $(I_n)$ is a decreasing sequence of intervals included in $(d,1]$ with $d>0$ satisfying \eqref{Gou}, then $B_n = \{ \theta^n \in I_n \}$ is strongly Borel-Cantelli as shown in Kim (2007, Prop. 4.1). 

We consider here two particular cases:
\begin{itemize}
 \item Consider  $I_n = [0,a_n]$ with $(a_n)_{n > 0}$ a decreasing sequence of real numbers in $]0,1]$ converging to $0$. Set $B_n = \{\theta^n \in I_n\} $. Using the same arguments as in Proposition 4.2 in Kim (2007), one can prove that, if $\sum_{n>0} a_n < \infty$, then $\mu ( \limsup_{n \rightarrow \infty} B_n) = 0$. Conversely, if $\sum_{n>0} a_n = \infty$, which is exactly condition \eqref{Gou}, then $(B_n)_{n \geq 1}$ is  a  Borel-Cantelli sequence. 
 
Now, to apply Theorem \ref{thmalphafaible} (and its Remark \ref{rmkthmalphafaible}), we first note that it has been proved in \cite{DGM}, that the $\tilde \alpha_{1, \infty} (n)$ coefficients of weak dependence associated with $(\theta^n)_{n \geq 0}$, viewed as a random sequence defined on $([0,1], \mu ) $, satisfy $  \kappa_1 n^{- ( 1- \gamma)/\gamma } \leq  \tilde \alpha_{1, \infty} (n) \leq  \kappa_2 n^{- ( 1- \gamma)/\gamma }$ for any $n \geq 1$ and some positive constants  $\kappa_1$ and $\kappa_2$.  Hence, in that case, Theorem \ref{thmalphafaible} gives the same condition \eqref{Gou} for the Borel-Cantelli property, up to the mild additional assumption  $n^{1/\gamma} a_n \rightarrow \infty$. This  shows that the approach based on the $\tilde \alpha_{1, \infty} (n)$  dependence coefficients provides optimal results in this case. Now, if  $na_n \rightarrow \infty$, then  $(B_n)_{n \geq 1}$ is a $L^1$ Borel-Cantelli sequence.  Finally, if 
 $\sum_{n \geq 1} n^{-1}  (n  a_n )^{\gamma -1}  < \infty$,  then  $(B_n)_{n \geq 1}$ is a strongly Borel-Cantelli sequence. 
 
\item Let  now $(a_n)_{n \geq 0}$ and  $(b_n)_{n \geq 0}$ be two  sequences of real numbers in $[0,1]$  such that $a_0 >0 $ and $b_{n+1} = b_n +a_n$ mod $1$. Define, for any $n \in \BBN$, 
 $I_{n+1} = [b_n,b_{n+1}]$  if $b_{n} < b_{n+1} $ and   $I_{n+1} =  [b_n, 1] \cup [0,b_{n+1}]$ if $b_{n+1} < b_{n} $.  It follows that   $(I_n)_{n \geq 1}$ is a sequence of consecutive intervals on the torus ${\mathbb R}/{\mathbb Z}$. Assume that  $\sum_{n \in \BBN} a_n = \infty$ (which is exactly \eqref{Gou}).  Since $\mu (I_{n+1}) \geq C a_n$, the divergence of the series implies that $\sum_{n>0} \mu (I_n) = \infty$.  Applying Theorem  \ref{BCbetafaible}  {\rm (iii)}, it follows that for any $\gamma < 1/2$, 
 $(B_n)_{n \geq 1}$ is a strongly Borel-Cantelli sequence.  Now if $\gamma = 1/2$, applying Theorem  \ref{BCbetafaible} {\rm (ii)} and {\rm (iii)} with $p=1$, 
 we get that  $(B_n)_{n \geq 1}$ is  a $L^1$ Borel-Cantelli sequence as soon as $(\sum_{k=1}^n a_k)/\log(n) \rightarrow \infty$,
 and a strongly Borel-Cantelli sequence as soon as $(\sum_{k=1}^n a_k)/(\log(n))^{2+ \varepsilon} \rightarrow \infty$ for some $\varepsilon>0$. 
 If $\gamma >1/2$, 
 we get that  $(B_n)_{n \geq 1}$ is  a $L^1$ Borel-Cantelli sequence as soon as $(\sum_{k=1}^n a_k)/n^{(2\gamma-1)/\gamma}\rightarrow \infty$, 
 and a strongly Borel-Cantelli sequence as soon as $(\sum_{k=1}^n a_k)/(n^{(2\gamma-1)/\gamma} (\log(n))^{1+ \varepsilon})\rightarrow \infty$
 for some $\varepsilon>0$.

\end{itemize}
\end{Example}

\begin{Example}\label{MC}
Let $(\varepsilon_i)_{ i \in \BBZ}$ be a sequence of iid random variables with values in $\BBR$, 
such that ${\mathbb E}(\log(1 + |\varepsilon_0|)) < \infty$. 
We consider here the stationary process
\beq \label{defAR1nonirreducible}
X_k= \sum_{i\geq 0} 2^{-i}  \varepsilon_{k-i} \,  ,
\eeq
which is defined almost surely (this is a  consequence of the three series theorem). The process $(X_k)_{k \geq 0}$ is a
 Markov chain, since  $X_{n+1}=\frac 12 X_n +  \varepsilon_{n+1}$.  However this  chain fails to be irreducible when 
 the innovations are with values in $\BBZ$. 
Hence the results of Sections \ref{sectionbetaalpha} and \ref{HarrisMC} cannot be applied in general. Nevertheless, under some mild additional conditions, the coefficients $\tilde \beta_{1,1}(n)$ of this chain converge to $0$ as shown by the lemma below.
\begin{Lemma}\label{easy}
Let $\mu$ be the law of $X_0$. Assume that $\mu$ has a bounded density.  If 
\begin{equation}\label{weakp}
 \textstyle  \sup_{t>0}   t^p {\mathbb P} (\log(1 + |\varepsilon_0|)>t)  < \infty  \quad \text{for some $p>1$}\, ,
\end{equation}
then $\tilde \beta_{1,1}(n) ={\cal O} ( n^{-(p-1)/2})$. 
\end{Lemma}
\begin{Remark}
The assumption that $\mu$ has a bounded density can be verified in many cases. 
For instance, it is satisfied if $\varepsilon_i= \xi_i + \eta_i$ where $(\xi_i)$ and $(\eta_i)$ are  two independent sequences of iid random variables, and 
$\xi_0$ has the Bernoulli$(1/2)$ distribution. 
Indeed, in that case, $X_0 = U_0 + Z_0$ with $U_0 =  \sum_{i=0}^\infty 2^{-i} \xi_{-i}$ and 
$Z_0 =  \sum_{i=0}^\infty  2^{-i} \eta_{-i}$. Since $U_0$ is uniformly distributed over $[0,2]$, it follows that the density of $\mu$ 
 is uniformly bounded by $1/2$.  
\end{Remark}

Since $(X_k)_{k \in {\mathbb Z}}$ is a stationary Markov chain, $ \tilde \alpha_{\infty, 1} (n) \leq \tilde \beta_{1,1}(n)$. Hence, under the assumptions  of Lemma \ref{easy}, we also have that $\tilde  \alpha_{\infty, 1}(n) ={\cal O} ( n^{-(p-1)/2})$. Let then $B_n = \{X_n \in I_n\}$.
As a consequence, we infer from Lemma \ref{easy},  Theorems \ref{BCbetafaible} and  \ref{thmalphafaible}   that

\begin{itemize}
\item 
If   $\mu  ( \limsup_n I_n ) >0$,  $\mu$ has a bounded density and \eqref{weakp} holds for some $p>3$, then $(B_n)_{n\geq 1}$  is a Borel-Cantelli 
sequence.
\item 
If $\mu$ has a bounded density,  \eqref{weakp} holds,   $\sum_{n\geq 1} \bigl( \mu ( I_n)  \bigr)^{(p+1)/(p-1)} = \infty$, $ (  \mu ( I_n) )_{n \geq 1}$ is non-increasing, and  $ \lim_{n } n^{(p-1)/2} \mu ( I_n)  = \infty$, 
 then
$(B_n)_{n\geq 1}$  is a Borel-Cantelli sequence.
\end{itemize}
\end{Example}

\begin{Example}\label{circle}
We consider the symmetric random walk on the circle, whose Markov 
 kernel is defined by 
\begin{equation}
  K f(x) = \frac{1}{2} \bigl( f (x+a) + f(x-a) \bigr)
\end{equation}
on the torus ${\mathbb R}/{\mathbb Z}$ with
 $a$ irrational in $[0,1]$. The Lebesgue-Haar  measure $\lambda$  is the unique
probability which is invariant by $K$. Let $(X_i)_{i\in \BBN}$ be
the stationary Markov chain with transition kernel $K$  and
invariant distribution $\lambda$.
We assume that $a$ is badly approximable in the weak sense meaning that, for any positive $\epsilon$, there exists some
positive constant $c$ such that
\begin{equation} \label{badly}
d(ka, {\mathbb Z} ) \geq  c |k|^{- 1 - \epsilon} \ \text{ for any }  k>0 .
\end{equation}
From Roth's theorem the algebraic numbers are badly approximable in the weak sense (see for instance   Schmidt
\cite{S}).  Note also that the set of numbers in $[0,1]$ satisfying \eqref{badly} has Lebesgue measure 1.
For this chain, we will obtain the bound below on the coefficients $\tilde \beta_{1,1}(n)$.

\begin{Lemma}\label{notsoeasy}
Let $a$ be badly approximable in the weak sense, and let $(X_i)_{i\in \BBN}$ be
the stationary Markov chain with transition kernel $K$  and
invariant distribution $\lambda$. Then, for any $b$ in $(0,1/2)$,  $\tilde \beta_{1,1}(n) = {\cal O} (n^{-b})$. 
\end{Lemma}

Since $(X_k)_{k \in {\mathbb Z}}$ is a stationary Markov chain, $ \tilde \alpha_{\infty, 1} (n) \leq \tilde \beta_{1,1}(n)$. Hence, under the assumptions  of Lemma \ref{notsoeasy},  $\tilde  \alpha_{\infty, 1}(n) ={\cal O} ( n^{-b})$ for any $b$ in $(0,1/2)$. 
As a consequence, we infer from Lemma \ref{notsoeasy},  Theorems \ref{BCbetafaible} and  \ref{thmalphafaible}   the corollary below
on the symmetric random walk on the circle with linear drift. 

\begin{Corollary} \label{RWwithdrift} Let $t$ be a real in $[0,1[$. Set $Y_k = X_k - kt$.  For any positive integer $n$, let 
$I_n = [0, n^{-\delta} ]$.  Set $B_n = \{ Y_n \in I_n \}$.  If $\delta < 1/3$,   $(B_n)_{n\geq 1}$  is a strongly Borel-Cantelli sequence for 
any $t$ in $[0,1[$.  Now, if $t$ is badly approximable in the strong sense, which means that (\ref{badly}) holds with $\epsilon = 0$,  
$(B_n)_{n\geq 1}$  is a strongly Borel-Cantelli sequence for any $\delta < 1/2$. 
\end{Corollary}

\end{Example}

\section{Harris recurrent Markov chains}  \label{HarrisMC}  

\setcounter{equation}{0}

In this section, we are interested in the Borel-Cantelli lemma for  irreducible and positively recurrent Markov chains.
Let $E$ be a Polish space and $\cal B$ be its Borel $\sigma$-field. 
Let $P$ be a  stochastic kernel. We assume that  there exists  a measurable function  $s$ with values in $[0,1]$ 
and a probability measure $\nu$ such that $\nu (s)>0$ and
\beq
\label{minorkernel}
P (x, A ) \geq s(x) \nu (A)
\ \text{ for any }\
(x,A) \in E \times {\cal B}.
\eeq
Then the chain is aperiodic and irreducible. Let us then define the sub-stochastic kernel $Q$ by 
\beq
\label{defQ}
Q(x,A) = P(x,A) - s(x) \nu (A)
\ \text{ for any }\
(x,A) \in E \times {\cal B}.
\eeq
Throughout this section, we assume furthermore that 
\beq
\label{Pitmanfinitemass}
 \sum_{n\geq 0} \nu Q^n  (1) <  \infty  . 
\eeq
Then the probability measure 
\beq
\label{invariantprobability}
\mu =   \Bigl( \sum_{n\geq 0} \nu Q^n  (1) \Bigr)^{-1} \sum_{n\geq 0}  \nu Q^n
\eeq
is the unique invariant probability measure under $P$. Furthermore the stationary Markov chain $(X_i)_{i \in \mathbb N}$ with kernel $P$ is positively recurrent 
(see Rio (2017), Chapter 9 for more details) and $\beta$-mixing according to Corollary 6.7 (ii) in Nummelin (1984).  Thus a direct  application of Theorem  \ref{BCbeta} (i) gives the following result.

\begin{Theorem} \label{BCHarrisbeta} Let $(A_k)_{k>0}$ be a sequence of Borel
subsets of $E$ such that $\mu ( \limsup_{n } A_n) >0$. 
Then  $\sum_{k>0}  \BBI_{A_k} (X_k)  = \infty $  a.s.   
\end{Theorem}

Obviously the result above does not apply in the case where the events are  nested and $\lim_n \mu (A_n) =0$. However in this case, the regeneration technique can be applied to prove the following result. 

\begin{Theorem} \label{BCHarris} Let $(A_k)_{k>0}$ be a sequence of Borel
subsets of $E$ such that $\sum_{k>0}  \nu (A_k) = \infty$ and $A_{k+1} \subset A_k$ for any positive $k$. 
Then  $\sum_{k>0}  \BBI_{A_k} (X_k)  = \infty $  a.s.   
\end{Theorem}

Suppose now that $\mu ( \limsup_{n } A_n) =0$ and that the events $(A_n)_{n \geq 1}$ are not necessarily nested. Then applying Corollary \ref{coralphadirect} and using Proposition 9.7 in Rio (2017) applied to arithmetic rates of mixing (see Rio (2017) page 164 and page 165 lines 8-11), we derive the following result:

\begin{Theorem} \label{BCHarrisalpha}  Let $T_0$ be the first renewal time of the extended Markov chain (see   \eqref{defTk} for the exact definition).  Assume that there exists $a >1$ such that 
$\BBP_{\mu} ( T_0 > n) \leq C n^{-a} $ for $n \geq 1$.  Suppose furthermore that  $(A_k)_{k>0}$ is a sequence of Borel
subsets of $E$ such that  $\sum_{n \geq 1}(  \mu (A_n) )^{(a+1)/a}= \infty$, $\lim_n n^a \mu (A_n)  = \infty$ and $(\mu (A_n) )_{n \geq 1}$ is non-increasing.  Then  $\sum_{k>0}  \BBI_{A_k} (X_k)  = \infty $  a.s.  
\end{Theorem}

If the stochastic kernel  $Q_1 ( x , . )$ defined in (\ref{defQ1}) is equal to $\delta_x$, then Theorem \ref{BCHarris} cannot be
further improved, as shown in Theorem \ref{BCHarrisconverse} below

\begin{Theorem} \label{BCHarrisconverse} Let $E$ be a Polish space. Let $\nu$ be a probability measure on 
$E$ and $s$ be a measurable function with 
values in $]0,1]$ such that $\nu (s) >0$. Suppose furthermore that 
\beq\label{Condposrec}
\int_E   {1\over s(x) }  d\nu (x) < \infty .
\eeq
Let  
\beq\label{ExampleDMR} 
P (x, . ) = s(x) \nu + (1- s(x) ) \delta_x . 
\eeq
Then $P$ is irreducible, aperiodic and positively recurrent. 
Let $(X_i)_{i \in \BBN}$ denote the strictly stationary Markov chain with kernel $P$ and $(A_k)_{k>0}$ be a sequence of Borel
subsets of $E$ such that $\sum_{k>0}  \nu (A_k) < \infty$ and $A_{k+1} \subset A_k$ for any positive $k$. 
Then  $\sum_{k>0}  \BBI_{A_k} (X_k)  < \infty $  a.s.   
\end{Theorem}

\begin{Remark} \label{Rioexample} Let us compare Theorems  \ref{BCHarris}  and  \ref{BCHarrisalpha} when $P$ is the Markov kernel defined by \eqref{ExampleDMR}  with
$E=[0,1]$,  $s(x) =x$ and  $\nu = (a+1)x^{a} \lambda$ with $a>0$ (here $\lambda$ is the Lebesgue measure on $[0,1]$).  For this example, 
 $\mu = ax^{a-1} \lambda$ and 
$\BBP_{\mu} ( T_0 > n)  \sim a \Gamma (a) n^{-a} $. Furthermore, from Lemma 2, page 75 in  Doukhan, Massart and Rio (1994), 
if $(\beta_n)_{n>0}$ denotes the sequence of $\beta$-mixing coefficients of the stationary Markov chain with kernel $P$, then 
\[ 
a \Gamma (a) \leq \liminf_n n^a \beta_n \leq \limsup_n n^a \beta_n \leq 3 a \Gamma (a) 2^a . 
\]
Now, for any $k \geq 1$, let  $A_k=I_k = ]a^{1/a}_k,b^{1/a}_k] $.
\begin{itemize}
\item Assume that  $I_{k+1} \subset I_k$, which means that $(a_k)$ is non-decreasing and $(b_k)$ is non-increasing.  
Then Theorem \ref{BCHarris} applies  if $\sum_{k>0} (b_k^{(a+1)/a} - a_k^{(a+1)/a} ) = \infty$ whereas Theorem \ref{BCHarrisalpha} applies if 
$\lim_n  n^a (b_n -a_n) = \infty$ and $\sum_{k>0} (b_k - a_k)^{(a+1)/a}  = \infty$.  Note that the first condition is always weaker than the second one. Note also that,  if 
$\lim_{k } a_k>0$, the first condition is equivalent to $\sum_{k >0} (b_k-a_k)= \infty$, which is then  strictly weaker than $\sum_{k>0} (b_k - a_k)^{(a+1)/a}  = \infty$.
Since $(b_k-a_k)= \mu(I_k)={\mathbb P}(X_k \in I_k)$,  the condition $\sum_{k>0} (b_k-a_k) = \infty$ is the best possible for the 
Borel-Cantelli property (this is due to the direct part of the Borel-Cantelli lemma).

\item Assume now that $a_k \equiv 0$ and $(b_k)$ is non-increasing. In that case, $\nu (I_k) = (  \mu (I_k) )^{(a+1)/a}$, for any $k \geq 1$. According to Theorem  \ref{BCHarrisconverse}, it follows that $\sum_{n \geq 1}(  \mu (I_n) )^{(a+1)/a}= \infty$ is a necessary condition to get the Borel-Cantelli property. 

\item Assume now that $I_k = ]a^{1/a}_k,(2a_k)^{1/a}] \subset [0,1]$ with $(a_k)_k \downarrow 0$. Since $I_{k+1} \not \subset I_k$ in this case, Theorem \ref{BCHarris} does not apply whereas the conditions of Theorem \ref{BCHarrisalpha} hold 
provided that $\lim_n n^a a_n = \infty$ and $\sum_{k>0} a_k^{(a+1)/a}  = \infty$. \end{itemize} 
\end{Remark}

\section{Proofs}

\setcounter{equation}{0}

\subsection{Proofs of the results of Section \ref{sectioncriteriageneral}} 

\subsubsection{Proof of Proposition \ref{lmacaracterisationBC}.}    We start by showing that $2. \Rightarrow 1.$ Let ${\Gamma} = \limsup_k \Gamma_k$. It suffices to prove that $ \BBP ({\Gamma} ) \geq \delta$. 
Note first that 
\[
\left\Vert \frac{   \sum_{k=1}^n   \BBI_{\Gamma_k} } {\sum_{k=1}^n \BBP  (  \Gamma_k)}  \right \Vert_{1} = 1 
\ \text{ and }\ 
\limsup_n  \left \Vert \frac{   \BBI_{\Gamma}  \sum_{k=1}^n   \BBI_{\Gamma_k} } {\sum_{k=1}^n \BBP  (  \Gamma_k)}  \right \Vert_{1}
 \leq \delta^{-1} \BBP  (  \Gamma) \, ,
\]
by \eqref{critsequenceGamma}. Hence it is enough to prove that 
\[
\lim_n  \left \Vert   \BBI_{\Gamma^c}  \frac{ \sum_{k=1}^n   \BBI_{\Gamma_k}  }{ \sum_{k=1}^n \BBP  (  \Gamma_k)}  \right \Vert_1 = 0 \, .
\]
This follows directly from \eqref{critsequenceGamma} and the fact that, by  definition of the  $\limsup$ and since $\sum_{k >0}  \BBP (\Gamma_k) = + \infty$, 
\[
\lim_n  \BBI_{\Gamma^c } \frac{ \sum_{k=1}^n   \BBI_{\Gamma_k} } {\sum_{k=1}^n \BBP  (  \Gamma_k)}   = 0  \mbox{ $ \BBP$-a.s.} 
\]
We prove now that $1. \Rightarrow 2.$    Proceeding by induction on $k$ one can construct  an
increasing sequence $(n_k)_{k\geq 0}$ of integers such that $n_0=1$  and 
\beq
\label{propertysequencenkbis}  
\BBP  \Bigl(   \bigcup_{j = n_{k-1} } ^ {n_k -1}  A_j \Bigr)  \geq \delta  ( 1 - 2^{-k})  \ \text{ for any }  k>0 . 
\eeq
Define now the sequence $(\Gamma_j)_{j > 0}$ of Borel sets by 
\[
\Gamma_{n_k} = A_{n_k}  \ \text{ and }  \Gamma_j =  A_j \setminus \Bigl( \bigcup_{i= n_k}^{j-1}  A_i \Bigr) \ \text{ for any } j \in ]n_k , n_{k+1} [ , 
\ \text{ for any }  k\geq 0 \, .
\]
 From the definition of $(\Gamma_j)_{j> 0}$
\[
\sum_{i=n_k}^{n_{k+1} - 1}  \BBI_{\Gamma_i }  = \BBI_{\bigl( \bigcup_{ \atop i \in [ n_k , n_{k+1} [ }  \Gamma_i  \bigr) }=  \BBI_{\bigl( \bigcup_{ \atop i \in [ n_k , n_{k+1} [ }  A_i  \bigr) }  \leq 1 \ \text{ for any }  k\geq 0 \,  . 
\]
Consequently, for any $j \geq 0$ and  any $n$ in $[n_j , n_{j+1} [$, 
\[
\sum_{i=1}^n \BBI_{\Gamma_i}  \leq \sum_{k=0}^j   \Bigl( \sum_{i=n_k}^{n_{k+1} - 1}  \BBI_{\Gamma_i } \Bigr) \leq j+1  \, . 
\]
Furthermore, from  (\ref{propertysequencenkbis}),
\[
\sum_{i=1}^n \BBP (\Gamma_i)  \geq \sum_{k=1}^j \BBP  \Bigl(   \bigcup_{i= n_{k-1} }^ {n_k -1}  A_i \Bigr)  \geq (j-1)\delta  
\]
for any $j \geq 1$ and  any $n$ in $[n_j , n_{j+1} [$. Hence, if $G_n =  \big (  \sum_{i=1}^n \BBP ( \Gamma_i) \big )^{-1} \sum_{i=1}^n \BBI_{\Gamma_i}$, then 
$G_n \leq (j+1)/((j-1)\delta)$ for $n$ in $[n_j , n_{j+1} [$, which ensures that $ \limsup_{n}G_n \leq 1/\delta$. 
\par\ssk
We now prove the second part of Proposition \ref{lmacaracterisationBC}.  Suppose that 
there exists a triangular sequence of events $(A_{k,n})_{1 \leq k \leq n}$ with $A_{k,n} \subset A_k$, such that ${\tilde E}_n
 = \sum_{k=1}^n \BBP (A_{k,n} ) \rightarrow \infty$ and that the sequence $(Z_n)_{n\geq 1}$ defined by 
  $Z_n =  {\tilde E}_n^{-1}  \sum_{k=1}^n \BBI_{A_{k,n}}  $ is uniformly integrable. Set $C_N = \bigcup_{ k > N }  A_k$. For any $n>N$, 
$$
\BBE ( Z_n)  =  \BBE \bigl( Z_n \BBI_{ C_N^c } \bigr)  +  \BBE \bigl( Z_n \BBI_{ C_N} \bigr) \leq 
( N / \tilde E_n ) +  \BBE \bigl( Z_n \BBI_{ C_N} \bigr) ,
$$
since $ \sum_{k=1}^n \BBI_{A_{k,n}}   \leq N$ on  $C_N^c$.  Using Lemma 2.1 (a) in Rio (2017), it follows that 
\[
1 = \BBE (Z_n) \leq (N/\tilde E_n)  +   \int_0^{1}  Q_{Z_n} (u)  Q_{ \BBI_{C_N}} (u)   du \\
 \leq (N/\tilde E_n)  +  \sup_{n>0}  \int_0^{\BBP ( C_N)}  Q_{Z_n} (u) du , 
\]
where $Q_Z$ denotes the cadlag inverse of the tail function  $t \mapsto \BBP ( Z>t)$. Hence,  
$$
1 = \lim_n  \BBE (Z_n)   \leq \sup_{n>0}  \int_0^{\BBP ( C_N)}  Q_{Z_n} (u) du .
$$
Now, if $\BBP  ( \limsup_k A_k)  = 0$, then $\lim_N \BBP ( C_N) = 0$. 
If furthermore $(Z_n)_{n>0}$ is  uniformly integrable, then, by Proposition \ref{propquantileuniformint}, the term on right hand
in the above inequality tends to $0$ as $N$ tends to $\infty$,  
which is a contradiction. The proof  of Proposition \ref{lmacaracterisationBC} is complete.  $\diamond$

\subsubsection{Proof of Corollary \ref{linkBCL1BC}.}   The fact that 2. implies 1. is immediate. Now, if 1. holds true,  then, by 
Proposition \ref{lmacaracterisationBC}, there exists a sequence $(\Gamma_k)_{k>0}$ of events 
 such that $\Gamma_k \subset A_k$, $\sum_{k > 0} \BBP  (  \Gamma_k) = + \infty$ and \eqref{critsequenceGamma} 
holds with $\delta=1$. 
Since $\Vert \sum_{k=1}^n   \BBI_{\Gamma_k} / \sum_{k=1}^n \BBP  (  \Gamma_k) \Vert_1 = 1$, it follows that 
\beq  \label{limnormsup}
\lim_n  \left \Vert \frac{ \sum_{k=1}^n   \BBI_{\Gamma_k} } {\sum_{k=1}^n \BBP  (  \Gamma_k)}  \right \Vert_{\infty} =  1  \, .
\eeq
Now 
\[
\left \Vert \frac{ \sum_{k=1}^n   \BBI_{\Gamma_k} } {\sum_{k=1}^n \BBP  (  \Gamma_k)}  -1 \right \Vert_1  = 
2 \left \Vert \left( \frac{ \sum_{k=1}^n   \BBI_{\Gamma_k} } {\sum_{k=1}^n \BBP  (  \Gamma_k)}  -1 \right)_+ \right \Vert_1  
\leq 2  \left (  \left \Vert \frac{ \sum_{k=1}^n   \BBI_{\Gamma_k} } {\sum_{k=1}^n \BBP  (  \Gamma_k)}  \right \Vert_{\infty}  - 1 \right )_+ ,
\]
which, together with \eqref{limnormsup}, implies that the above sequence $(\Gamma_k)_{k>0}$ is a $L^1$ Borel-Cantelli sequence.
Hence Corollary \ref{linkBCL1BC} holds.  $\diamond$

\subsubsection{Proof of Proposition \ref{lmacaracterisationBCint}.}    The fact that  $2. \Rightarrow 1.$ follows immediately from 
Proposition \ref{lmacaracterisationBC}. We now prove the direct part. Proceeding by induction on $k$ one can construct  an
increasing sequence $(n_k)_{k\geq 0}$ of integers such that $n_0=1$  and 
\beq
\label{propertysequencenkbisint}  
\mu \Bigl(   \bigcup_{j = n_{k-1} } ^ {n_k -1}  I_j \Bigr)  \geq \delta  ( 1 - 2^{-k})  \ \text{ for any }  k>0 . 
\eeq
Now, for any $k\geq 0$, we construct the intervals $\Gamma_j$ for $j$ in $[n_k, n_{k+1}[$. This will be done by using the lemma below. 

\begin{Lemma} \label{Coveringdisjointint}  Let $(J_k)_{k\in [1,m]}$ be a sequence of intervals of $\BBR$. Then there exists a sequence $(\Gamma_k)_{k\in [1,m]}$ of disjoint intervals such that $\bigcup_{k=1}^m \Gamma_k = \bigcup_{k=1}^m J_k$ and 
$\Gamma_k \subset J_k$ for any $k$ in $[1,m]$.  
\end{Lemma}
\noindent
{\bf Proof of Lemma \ref{Coveringdisjointint}.}  We prove the Lemma by induction on $m$. Clearly the result holds true for $m=1$. 
Assume now that Lemma \ref{Coveringdisjointint} holds true at range $m$. Let then $(J_k)_{k\in [1,m+1]}$ be a sequence of intervals.
By the induction hypothesis, there exists a sequence $(\Gamma_{k,m})_{1\leq k \leq m}$ of disjoint intervals such that 
$\bigcup_{k=1}^m \Gamma_{k,m} = \bigcup_{k=1}^m J_k$ and $\Gamma_{k,m} \subset J_k$ for any $k$ in $[1,m]$.  Now, at the range $m+1$, 
define now  the intervals $\Gamma_k$ for $k$ in $[1,m]$ by $\Gamma_k = \emptyset$ if $\Gamma_{k,m} \subset J_{m+1}$ and 
$\Gamma_k = \Gamma_{k,m}$  if $\Gamma_{k,m} \not\subset J_{m+1}$. Clearly these intervals are disjoint. Set 
\beq
\label{defGamma(m+1)}
\Gamma_{m+1} =  \bigcap_{k=1}^m \Bigl( \Gamma_k^c \cap J_{m+1} \Bigr) . 
\eeq
If $\Gamma_k = \emptyset$, then $\Gamma_k^c \cap J_{m+1} = J_{m+1}$. Otherwise, from the definition of $\Gamma_k$, 
$\Gamma_k$ is a nonempty interval and $\Gamma_k \not\subset J_{m+1}$, which implies that $\Gamma_k^c \cap J_{m+1} $
is an interval.  Hence $\Gamma_{m+1}$ is a finite intersection of intervals, which ensures that $\Gamma_{m+1}$ is an interval.
By \ref{defGamma(m+1)},   $\Gamma_{m+1}$ does not intersect $\Gamma_k$ for any $k$ in $[1,m]$.
Hence the so defined intervals $\Gamma_k$ are disjoint, $\Gamma_k \subset J_k$  for any $k$ in $[1,m+1]$. 
Finally
\beq
\label{unionGammak}
\bigcup_{k=1}^{m+1}  \Gamma_k  = J_{m+1}  \bigcup \Bigl( \bigcup_{k=1}^m \Gamma_k \Bigr) =  
J_{m+1}  \bigcup \Bigl( \bigcup_{k=1}^m \Gamma_{k,m} \Bigr) =  J_{m+1}  \bigcup \Bigl( \bigcup_{k=1}^m J_k \Bigr) 
\eeq
Hence, if Lemma \ref{Coveringdisjointint} holds true at range $m$, then Lemma \ref{Coveringdisjointint} holds true at range $m+1$,
which ends the proof of the lemma.   $\diamond$
\par\medskip\no

\smallskip

\noindent{\bf End of the proof of Proposition \ref{lmacaracterisationBCint}.} For any $k\geq 0$,  by Lemma \ref{Coveringdisjointint} applied to  
 $(I_j)_{j\in [n_k , n_{k+1}[}$, there exists a sequence $(\Gamma_j)_{j\in [n_k, n_{k+1} [ }$ of disjoint intervals such that 
\beq
\bigcup_{j\in [n_k , n_{k+1} [ }  \Gamma_j =  \bigcup_{j\in [n_k , n_{k+1} [ }  I_j  \ \text{ and } \Gamma_j \subset I_j \ \text{ for any }
j \in [n_k, n_{k+1} [ . 
\eeq
From now on the end of the proof is exactly the same as the end of the proof of the first part of Proposition \ref{lmacaracterisationBC}. $\diamond$

\subsubsection{Proof of Proposition \ref{PropCriteriaBC}.}   We start by proving Item (ii).  Let $f$ be the function  defined in Definition \ref{deffunctionf} and $X$ be any integrable real-valued  random variable. Then
\beq\label{comparisonL1Ef}  
\Vert X \Vert_1 \leq  \Vert X \BBI_{|X| \leq 1} \Vert_2 +  \Vert X \BBI_{|X|>1}  \Vert_1   \leq \sqrt{ 2 \BBE ( f(X) )}  +  2 \BBE ( f(X) ) .
\eeq
Consequently, if (\ref{CritL1BC}) holds, then $\lim_{n\rightarrow \infty} \Vert (S_n-E_n)/E_n \Vert_1 = 0$, which proves Item (ii). 
\par\smallskip
\noindent{\bf Proof of Item (i).}   Applying (\ref{comparisonL1Ef}),  we get that 
$\lim_{k\rightarrow \infty}  \Vert (\tilde S_{n_k} /\tilde E_{n_k} ) - 1 \Vert_1 = 0$. Hence, by the Markov inequality, 
$\lim_{k\rightarrow \infty}  \BBP  ( \tilde S_{n_k} \leq \tilde E_{n_k} /2 ) = 0$, which proves that $\tilde S_{n_k}$ converges to $\infty$
in probability as $k$ tends to $\infty$. 
Now  $g_{j,n_k} \leq \BBI_{ B_j}$ any $j$ in $[1, n_k]$. Therefrom $\tilde S_{n_k}  \leq S_{n_k}$ and consequently $S_{n_k}$ converges to $\infty$
in probability as $k$ tends to $\infty$. Since $(S_n)_n$ is a non-decreasing sequence of random variables, it implies immediately that 
$\lim_{n\rightarrow \infty}  S_n = +\infty$ almost surely, which completes the proof of Item (i).  
\par\smallskip
\noindent{\bf Proof of Item (iii).} For any non-negative real $x$, define $E : x \mapsto E(x) = \BBE (S_{[x]}) $. $E$ is a non-decreasing and cadlag function defined on ${\mathbb R}^+$ with values in ${\mathbb R}^+$. 
Let $E^{-1}$ be its generalized inverse on ${\mathbb R}^+$ defined by $E^{-1} (u) = \inf \{ x \in {\mathbb R}^+ \, : \, E(x) \geq u \}$. Hence 
\beq \label{proprieteinversedeE}
x \geq E^{-1} (u) \iff E(x) \geq u  .
\eeq 
Note that $E ([x]) = E_{[x]}$.  Let $\tau_n = \alpha^n$ for a fixed $\alpha >1$ and define 
\[
m_n = E^{-1} ( \tau_n) = \inf \{ k \geq 1 \, : \, E (k) \geq \tau_n \}\, .
\]
Hence $(m_n)_{n \geq 1}$ is a non-decreasing sequence of integers. Note also that there exists a positive integer $n_0 $ depending on $\alpha$ such that, for any $n \geq n_0$, $m_n < m_{n+1}$. Indeed, let assume that there exists $n \geq n_0$ such that $m_n =m_{n+1}$.  By definition $E(m_n -1) < \alpha^n$ and  $E(m_n) =E(m_{n+1}) \geq  \alpha^{n+1}$. This implies that 
\[
 \alpha^{n+1} \leq E(m_n -1) + \BBP (B_{m_n}) < \alpha^n + 1 \, . 
\]
Since $\alpha >1$, there exists an integer $n_0$ such that the above inequality fails to hold for any $n \geq n_0$. This contradicts the fact that there exists $n \geq n_0$ such that $m_n =m_{n+1}$. 
 Let us then show that 
\beq \label{CritstrongBCp1}
( S_{m_n} /  E_{m_n} ) \rightarrow 1 \text{ almost surely, as } \, n \rightarrow \infty \, .
\eeq
By the first part of the Borel-Cantelli lemma, \eqref{CritstrongBCp1} will hold provided that 
\beq \label{CritstrongBCp2}
\sum_{n \geq n_0} \BBE  \bigl(  f \bigl( (S_{m_n}-E_{m_n}) / E_{m_n}  \bigr) \bigr) <  \infty \, .
\eeq
Hence, setting, for any real $b >0$, 
\[
f^* (x /b ) := \sup_{1 \leq k \leq [x]} \BBE  \bigl(  f \bigl( (S_{k}-E_{k}) /b  \bigr) \bigr) \, , 
\]
to prove \eqref{CritstrongBCp2}, it suffices to show that
\beq \label{CritstrongBCp3}
\sum_{n \geq n_0}  f^* (m_n / E_{m_n} ) <  \infty \, .
\eeq
Write 
\begin{align*}
\sum_{n \geq n_0}  f^* (m_n / E_{m_n} )  & = \sum_{n \geq n_0}   \sum_{k=m_n+1}^{m_{n+1}} \BBP (B_k) f^* (m_n / E_{m_n} ) \Bigl (\sum_{k=m_n+1}^{m_{n+1}} \BBP (B_k)  \Bigr )^{-1} \\
& \leq  \sum_{n \geq n_0}   \sum_{k=m_n+1}^{m_{n+1}} \BBP (B_k) f^* (k  / E_{m_n} ) \Bigl (\sum_{k=m_n+1}^{m_{n+1}} \BBP (B_k)  \Bigr )^{-1}  \, .  
\end{align*}
Note now that, for any real $a \geq 1$,  $f(ax) \leq a^2 f(x)$. Therefore
\[
\sum_{n \geq n_0}  f^* (m_n / E_{m_n} )  
 \leq  \sum_{n \geq n_0}   \sum_{k=m_n+1}^{m_{n+1}} \BBP (B_k) (E_k/E_{m_n})^2 f^* (k  / E_{k} ) \Bigl (\sum_{\ell=m_n+1}^{m_{n+1}} \BBP (B_\ell)  \Bigr )^{-1}  \, .  
\]
Next, for any $k \leq m_{n+1}$, $E_k \leq E_{m_{n+1}} < \tau_{n+1} + \BBP (B_{m_{n+1}} ) \leq \alpha^{n+1} +1$, $E_{m_n} \geq \alpha^n$
and  
\[
\sum_{\ell=m_n+1}^{m_{n+1}} \BBP (B_\ell)  =   E_{m_{n+1}}  - E_{m_n} \geq \alpha^{n+1} - (  \alpha^n + \BBP (B_{m_{n}} ) )  \geq \alpha^n (\alpha -1) -1 \geq \alpha^n (\alpha -1)/2 \, ,  
\]
for any $n \geq n_1$.  Hence, for any $ n \geq n_{1}$ and any $1 \leq k \leq m_{n+1}$, 
\[
\sum_{\ell=m_n+1}^{m_{n+1}} \BBP (B_\ell)    \geq \frac{\alpha^{n+1} (\alpha -1)}{2 \alpha} \geq \frac{(\alpha^{n+1} +1) (\alpha -1)}{4 \alpha}  \geq \frac{E_k(\alpha -1)}{4 \alpha} \, .  
\]
So, overall, setting $n_2= \max (n_0,n_1)$,
\[
\sum_{n \geq n_2}  f^* (m_n / E_{m_n} )  
 \leq \frac{4 \alpha (\alpha +1)^2}{ \alpha -1}\sum_{n \geq n_2}   \sum_{k=m_n+1}^{m_{n+1}}  \frac{\BBP (B_k) }{E_k}f^* (k  / E_{k} )  =: C_{\alpha} \sum_{k \geq m_{n_2}}\frac{\BBP (B_k) }{E_k}f^* (k  / E_{k} )  \, ,  
\]
proving \eqref{CritstrongBCp3} (and subsequently \eqref{CritstrongBCp1}) under \eqref{CritstrongBC}.  The rest of the proof is quite usual but
we give it for completeness. Since $(S_n)_{n \geq 1}$ is a non-decreasing sequence as well as the normalizing sequence $(E_n)_{n \geq 1}$, if $1 < m_n \leq k \leq m_{n+1}$,
\[
\frac{E_{m_n}}{E_{m_{n+1}}}\frac{S_{m_{n}}}{ E_{m_n}} \leq \frac{S_k}{ E_k} \leq\frac{E_{m_{n+1}}}{E_{m_n}}\frac{S_{m_{n+1}}}{ E_{m_{n+1}}} \, .
\]
But, for any positive integer $k$, $\alpha^{k} \leq E_{m_{k}} < \alpha^{k} + \BBP(B_k)$. Therefore $E_{m_{n+1}}/E_{m_n} \rightarrow \alpha$, as $n \rightarrow \infty$.  Hence, by using \eqref{CritstrongBCp1}, almost surely, 
\[
(1/\alpha) \leq \liminf_{k \rightarrow \infty} (S_k / E_k ) \leq  \limsup_{k \rightarrow \infty} (S_k / E_k)  \leq \alpha \, .
\]
Taking the intersection of all such events for rationals $\alpha >1$, Item (iii) follows.   $\diamond$

\subsection{Proofs of the results of Section \ref{sectionbetaalpha}}

 \subsubsection{Proof of Theorem \ref{BCbeta} ($\beta$-mixing case)}   
 
Throughout this section,  $\beta_j = \beta_{\infty, 1} (j)$. Items (i) and (ii) will be derived from the proposition below.

\begin{Proposition} \label{Propbeta} With the notations of Theorem \ref{BCbeta}, let $(\Gamma_{k,n})_{1\leq k \leq n}$  be a double array of Borel sets in $E$. Set $\tilde E_n = \sum_{k=1}^n  \mu ( \Gamma_{k,n} )$ and  $G_n = \tilde E_n^{-1} \sum_{k=1}^n \BBI_{\Gamma_{k,n}}$. Suppose that 
$\tilde E_n >0$ for any positive $n$,  $\lim_{n\uparrow \infty} \tilde E_n = \infty$ and $(G_n)_{n>0}$ is a uniformly integrable sequence in $(E , {\cal B} (E) , \mu)$.  Let $B_{k,n} = \{X_k \in \Gamma_{k,n} \}$ and $\tilde S_n = \sum_{k=1}^n \BBI_{B_{k,n} }$. If $\, \lim_{n\uparrow \infty}   \beta_n = 0$, then 
\beq\label{PropbetaCVL1}   \lim_{n\rightarrow \infty}  \Vert  (\tilde S_n -\tilde E_n) / \tilde E_n \Vert_1 = 0 .   \eeq
\end{Proposition}

 \noindent
{\bf Proof of Proposition \ref{Propbeta}. }  
From (\ref{comparisonL1Ef}),  it is enough to prove that 
\beq\label{CritPropbeta}  
\lim_{n\rightarrow \infty} \BBE  \bigl( f_n ( \tilde S_n - \tilde E_n ) \bigr)= 0 , \ \text{where}\ f_n (x) = f (x/\tilde E_n) . 
\eeq
Now, by setting $ \tilde S_0= \tilde E_0=0$, we first write
\beq\label{decomp1}
f_n ( \tilde S_n - \tilde E_n )   =  \sum_{k=1}^n  \bigl(  f_n ( \tilde S_k - \tilde E_k) ) - f_n ( \tilde S_{k-1} - \tilde E_{k-1} ) \,  \bigr) . 
\eeq
Let then $T_0=0$ and, for $k>0$,
\beq \label{defxik}
T_k = \tilde S_k - \tilde E_k, \ \xi_k = T_k - T_{k-1} =   \BBI_{\Gamma_{k,n}} (X_k)  - \mu (\Gamma_{k,n})  . 
\eeq
With these notations, by the Taylor integral formula at order $1$, 
\begin{eqnarray*}
 f_n ( \tilde S_k - \tilde E_k) ) -   f_n ( \tilde S_{k-1} - \tilde E_{k-1} )  & = &  f_n ( T_k  ) - f_n ( T_{k-1}  ) 
\cr & = &  f'_n ( T_{k-1} )  \xi_k  +   \int_0^1  
\bigl( f'_n ( T_{k-1} + t \xi_k ) - f'_n ( T_{k-1}) \bigr) \xi_k \, dt  . 
\end{eqnarray*}
Now $f'_n (x) = \tilde E_n^{-1} f' ( x/\tilde E_n)$. Moreover, from the definition of $f$, $f'$ is $1$-Lipschitzian. Hence  
$$
\bigl( f'_n ( T_{k-1} + t \xi_k ) - f'_n ( T_{k-1}) \bigr) \xi_k    \leq   \tilde E_n^{-2}  \xi_k^2
$$
for any $t$ in $[0,1]$, which implies that 
\beq \label{step1beta}
f_n ( T_k  ) - f_n ( T_{k-1}  )   \leq    f'_n ( T_{k-1} )  \xi_k  +  \tilde E_n^{-2}   \xi_k^2   . 
\eeq
Now, using (\ref{decomp1}), (\ref{step1beta}), taking the expectation and noticing that $f'_n ( T_0 ) = f'_n (0) = 0$, we get that 
\beq\label{step2beta}
 \BBE  \bigl(  f_n ( \tilde S_n - \tilde E_n ) \bigr)  \leq \sum_{k=2}^n \BBE \bigl(  f'_n ( T_{k-1} ) \xi_k \bigr)  +  
\tilde E_n^{-2} \sum_{k=1}^n \mu (\Gamma_{k,n})  . 
\eeq
 Next, let $m\geq 2$ be a fixed integer.  For $n\geq m$,   
$$
 f'_n ( T_{k-1} ) \xi_k  = f'_n ( T_{(k -m)_+ } )  \xi_k  + \sum_{j=1}^{m-1} \bigl( f'_n ( T_{k-j} ) -  f'_n ( T_{(k-j-1)_+ } ) \bigr)  \xi_k  . 
$$
Taking the expectation in the above equality, we then  get that 
\beq\label{decomp2}
\BBE ( f'_n ( T_{k-1} ) \xi_k )   = \cov ( f'_n ( T_{(k -m)_+ } ) , \BBI_{B_{k,n} } )  + 
\sum_{j=1}^{m-1} \cov \bigl( f'_n ( T_{(k-j)_+} ) -  f'_n ( T_{(k-j-1)_+ } )  ,  \BBI_{B_{k,n} } )  . 
\eeq
In order to  bound up the terms appearing in (\ref{decomp2}), we will use Delyon's covariance inequality, which we now recall. We refer to Rio (2017, Theorem 1.4) for an  available reference with a proof. 

\begin{Lemma}  \label{Del90} - {\sl Delyon (1990)} - Let $\cal A$ and $\cal B$ be two $\sigma$-fields of
$\EP$. Then there exist  random variables $d_{\cal A}$ and $d_{\cal B}$ respectively  $\cal A$-measurable with values in 
$[0, \varphi ( {\cal A} , {\cal B} ) ]$ and  $\cal B$-measurable with values in 
$[0, \varphi ( {\cal B} , {\cal A} ) ]$ ,  satisfying 
$\BBE ( d_{\cal A} )  = \BBE ( d_{\cal B} ) =   \beta ({\cal A} , {\cal B} )$
and such that, for any  $(p,q)$ in $[1,\infty]^2$ with $(1/p)+ (1/q)= 1$ and any random vector
 $(X,Y)$  in $L^p ({\cal A}) \times L^q ({\cal B})$,
\beq \label{Delyon}
| \cov (X, Y) | \leq 2 \bigl( \BBE (d_{\cal A}  |X|^p ) \bigr)^{1/p}
\bigl( \BBE (d_{\cal B}  |Y|^q ) \bigr)^{1/q} , 
\eeq
where $\bigl( \BBE (d_{\cal A}  |X|^p ) \bigr)^{1/p} = \Vert X \Vert_\infty$ if $p= \infty$ and  
$\bigl( \BBE (d_{\cal B}  |Y|^q ) \bigr)^{1/q}  = \Vert Y \Vert_\infty$ if $q= \infty$. 
\end{Lemma}

We now bound up the first term in  the right-hand side of equality  (\ref{decomp2}). If $k\leq m$, then $T_{(k -m)_+} = 0$, whence
$$
\cov \bigl( f'_n ( T_{(k -m)_+ } ) , \BBI_{B_{k,n}}  \bigr) = 0 .
$$
Set
\beq \label{DefWkj}
W_{k,l}  =  \sum_{-k < i \leq -l} \bigl( \, \BBI_{\Gamma_{i+k,n} } (X_i)   - \mu ( \Gamma_{i+k,n} )  \bigr) \ \text{ for any }  l \leq  k . 
\eeq
If $k>m$, using  the stationarity of $\XZ$, we obtain that 
\beq \label{shiftcov}
\BBE \bigl( f'_n ( T_{(k -m)_+ } )  \xi_k  \bigr)  = \cov ( f'_n  ( W_{k,m}  )  , \BBI_{\Gamma_{k,n}} ( X_0) ) .
\eeq
Let us now apply Lemma \ref{Del90} with ${\cal A} = {\cal F}_{-m}$, ${\cal B} = \sigma (X_0)$,  $p=\infty$, $q=1$, 
$X = f'_n  ( W_{k,m} )$ and $Y = \BBI_{\Gamma_{k,n} } ( X_0)$:  
there exists some measurable function $\psi_m$ satisfying 
\beq \label{intpsi}
0 \leq \psi_m \leq  \varphi ( \sigma (X_0) , {\cal F}_{-m} )     \ \text{ and }\ \int_E \psi_m d\mu = \beta_m ,
\eeq
 such that, for any $k>m$, 
\beq \label{applicovDel1}
\cov ( f'_n ( T_{(k -m)_+ } ) , \BBI_{B_{k,n} } )   \leq 2  \Vert f'_n \Vert_\infty \int_E  \BBI_{\Gamma_{k,n}}  \psi_m d\mu . 
\eeq
Next $f'_n (x) = \tilde E_n^{-1} f' (x/\tilde E_n)$. Since $\Vert f' \Vert_\infty \leq 1$, it follows that $\Vert f'_n \Vert_\infty \leq \tilde E_n^{-1}$. 
Summing (\ref{applicovDel1}) on $k$ and using this bound, we finally get that 
\beq \label{term1}
\sum_{k=2}^n   \cov ( f'_n ( T_{(k -m)_+ } ) , \BBI_{B_{k,n} } )  \leq  2  \int_E  G_n  \psi_m d\mu ,
\eeq
where $G_n$ is defined in Proposition \ref{Propbeta}. 
\par
We now bound up the other terms in  the right-hand side of equality  (\ref{decomp2}).  If $j\geq k$, then 
$T_{(k-j)_+}  =   T_{(k-j-1)_+ } = 0$, which implies that 
$$
\cov \bigl( f'_n ( T_{(k-j)_+} ) -  f'_n ( T_{(k-j-1)_+ } )  ,  \BBI_{B_{k,n} } ) = 0 . 
$$
If $j<k$, using  the stationarity of $\XZ$, we obtain that 
\beq \label{shiftcovbis}
\cov \bigl( f'_n ( T_{(k-j)_+} ) -  f'_n ( T_{(k-j-1)_+ } )  ,  \BBI_{B_{k,n} } )  =  \cov ( f'_n  ( W_{k,j} )-f'_n ( W_{k,j+1}  )  , \BBI_{\Gamma_{k,n}} ( X_0) ) ,
\eeq
 where $W_{k,j}$ and $W_{k,j+1}$ are defined in (\ref{DefWkj}).  Applying Lemma \ref{Del90} with ${\cal A} = {\cal F}_{-j}$, 
${\cal B} = \sigma (X_0)$,  $p=q=2$,  $X =  f'_n  ( W_{k,j}) -f'_n ( W_{k,j+1}  ) $ and $Y = \BBI_{\Gamma_{k,n} } ( X_0)$,   we obtain that there exist
some $\sigma (X_0)$-measurable random variable $b_j$ and some ${\cal F}_{-j}$-measurable random variable $\eta_j$  
with values in $[0,1]$,  satisfying 
\beq \label{Ebjetaj}
\BBE ( b_j ) = \BBE ( \eta_j)   =  \beta_j
\eeq
and such that 
\beq  \label{applicovDel2}
\cov ( f'_n  ( W_{k,j} ) -f'_n ( W_{k,j+1}  )  , \BBI_{\Gamma_{k,n}} ( X_0) )  \leq 
2  \sqrt { \BBE (\eta_j | f'_n  ( W_{k,j} ) -f'_n ( W_{k,j+1}  )   |^2 ) 
 \BBE (b_j \BBI_{\Gamma_{k,n} }(X_0) ) }  \,  . 
\eeq
Next, from the definitions of $f_n$ and $f$, $f'_n (x) = \tilde E_n^{-1} f' (x/\tilde E_n)$ and $f'$ is $1$-Lipschitzian. Consequently 
$$
| f'_n  ( W_{k,j} ) -f'_n ( W_{k,j+1}  )   |   \leq \tilde E_n^{-2}  | W_{k,j}  - W_{k,j+1} |  =  
\tilde E_n^{-2} | \BBI_{\Gamma_{k-j,n} } (X_{-j} ) - \mu ( \Gamma_{k-j,n} ) |  , 
$$
which implies that 
$$
\BBE (\eta_j | f'_n  ( W_{k,j} ) -f'_n ( W_{k,j+1}  )   |^2 ) \leq \tilde E_n^{-4} 
\BBE ( b'_j  | \BBI_{\Gamma_{k-j,n} } (X_{-j} ) - \mu ( \Gamma_{k-j,n} ) |^2 ) ,
$$
with $b'_j = \BBE ( \eta_j  \mid \sigma (X_{-j} )  ) $. Combining the above inequality, (\ref{applicovDel2}) and the elementary
inequality $2 \sqrt{ab} \leq a+b$, we infer that 
\beq  \label{applicovDel2bisMarch}
\tilde E_n^2 \cov ( f'_n  ( W_{k,j})  -f'_n ( W_{k,j+1}  )  , \BBI_{\Gamma_{k,n}} ( X_0) )  \leq \BBE (b_j \BBI_{\Gamma_{k,n} }(X_0)  + b'_j  | \BBI_{\Gamma_{k-j,n} } (X_{-j} ) - \mu ( \Gamma_{k-j,n} ) |^2 ) .
\eeq
Recall now that $b_j$ is $\sigma (X_0)$-measurable and $b'_j$ is $\sigma (X_{-j} )$-measurable. Hence there exists  Borelian functions
$\varphi_{j,0}$ and $\varphi_{j,1}$ with values in $[0,1]$ such that $b_j = \varphi_{j,0} (X_0)$ and $b'_j = \varphi_{j,1} (X_{-j})$. 
Using now the stationarity of $\XZ$, we get
$$
\BBE (b_j \BBI_{\Gamma_{k,n} }(X_0) + b'_j  | \BBI_{\Gamma_{k-j,n} } (X_{-j} ) - \mu ( \Gamma_{k-j,n} ) |^2 )   = 
\int_E \bigl(  \varphi_{j,0} \BBI_{\Gamma_{k,n}}  + \varphi_{j,1} | \BBI_{\Gamma_{k-j,n} }  - \mu ( \Gamma_{k-j,n} ) |^2 \bigr) d\mu .
$$
Next, applying the elementary inequality 
$$
|\BBI_{\Gamma_{k-j,n} }  - \mu ( \Gamma_{k-j,n} ) |^2\leq \BBI_{\Gamma_{k-j,n} }  + \mu ( \Gamma_{k-j,n} ) ,
$$
noticing that $\int_E \varphi_{j,1} d\mu = \beta_j$ and putting together (\ref{shiftcovbis}),  (\ref{applicovDel2bisMarch}) 
and the above inequalities, we get
\beq
\label{applicovDel3}
\tilde E_n^2 \cov \bigl( f'_n ( T_{(k-j)_+} ) -  f'_n ( T_{(k-j-1)_+ } )  ,  \BBI_{B_{k,n} } ) \leq 
\beta_j \mu ( \Gamma_{k-j,n} ) +   \int_E  \bigl( \varphi_{j,0} \BBI_{\Gamma_{k,n}}  + \varphi_{j,1}  \BBI_{\Gamma_{k-j,n} } \bigr) d\mu,
\eeq
for some Borelian functions $\varphi_{j,0}$ and $\varphi_{j,1}$ with values in $[0,1]$ satisfying 
\beq
\label{Intvarphij01}
\int_E \varphi_{j,0} d\mu = \int_E  \varphi_{j,1} d\mu = \beta_j. 
\eeq
Finally, summing (\ref{applicovDel3}) on $j$ and $k$,  using (\ref{step2beta}), (\ref{decomp2}) and (\ref{term1}), we obtain
\beq  \label{Laststepbeta}
 \BBE  \bigl(  f_n ( \tilde S_n - \tilde E_n ) \bigr)  \leq \tilde E_n^{-1}  \Bigl( 1 + \sum_{j=1}^{m-1} \beta_j \Bigr) +
2 \int_E  G_n \Bigl(  \psi_m +  \tilde E_n^{-1} \sum_{j=1}^{m-1} \psi'_j \Bigr)  d\mu , 
\eeq
where 
\beq \label{defpsi'}
\psi'_j = (\varphi_{j,0}  + \varphi_{j,1} )/2 . 
\eeq 
Let  $\lambda$ denote the Lebesgue measure on $[0,1]$ and let $Q_{G_n}$ be the cadlag inverse function of the the tail function of $G_n$. Then, 
by Lemma 2.1 (a) in Rio (2017) 
applied to the functions  $G_n$ and $\BBI_{u\leq \psi_m}$, 
\beq \label{quantile1}
 \int_E  G_n  \psi_m d\mu  = \int \!\!\!\int_{E\times [0,1]}  G_n \BBI_{u\leq \psi_m}  d \mu\otimes \lambda  
\leq \int_0^{\beta_m} Q_{G_n} (s) ds . 
\eeq
In a similar way 
\beq \label{quantile2}
 \int_E  G_n  \psi'_j d\mu  
\leq \int_0^{\beta_j} Q_{G_n} (s) ds . 
\eeq
 Putting the two above inequalities in (\ref{Laststepbeta}), we get: 
\beq  \label{Laststepbetabis}
 \BBE  \bigl(  f_n ( \tilde S_n - \tilde E_n ) \bigr)  \leq \tilde E_n^{-1}  \Bigl( 1 + \sum_{j=1}^{m-1} 
\int_0^{\beta_j}  ( 1 + 2 Q_{G_n} (s) ) ds  \Bigr) +
2 \int_0^{\beta_m}  Q_{G_n} (s) ds .  
\eeq
\par
We now complete the proof of Proposition \ref{Propbeta}.  Since 
$\int_0^1 Q_{G_n} (s) ds = \int_E G_n d\mu = 1$,
the above inequality ensures that
\beq  \label{EndproofPropbeta2}
 \BBE  \bigl(  f_n ( \tilde S_n - \tilde E_n ) \bigr)  \leq \tilde E_n^{-1} (3m-2)   + 2 \int_0^{\beta_m} Q_{G_n} (s) ds .   
\eeq
It follows that 
\beq \label{EndproofPropbeta3}
\limsup_{n\rightarrow \infty } \BBE  \bigl( f \bigl( (\tilde S_n - \tilde E_n)/\tilde E_n \bigr) \bigr)  \leq    2  \limsup_{n \rightarrow \infty}
 \int_0^{\beta_m} Q_{G_n} (s) ds  
\eeq
for any  integer $m\geq 2$. 
Now $\lim_{m\uparrow \infty}   \beta_m = 0$. Consequently, if the sequence $(G_n)_{n>0}$ is uniformly integrable, 
then, by Proposition \ref{propquantileuniformint},  the term on right hand in the above inequality tends to $0$ as $m$ tends to $\infty$, which 
ends the proof of Proposition \ref{Propbeta}.  $\diamond$
\par\ssk

\smallskip

\noindent {\bf End of the proof of Theorem \ref{BCbeta}.}  Item (ii) follows immediately from Proposition  \ref{Propbeta} applied with $\Gamma_{k,n} = A_k$.  
To prove Item (i), we note that applying Proposition \ref{lmacaracterisationBC} with $(\Omega, {\mathcal T}, \BBP) = (X, {\mathcal B} (X) , \mu)$, there exists a sequence of events $(\Gamma_k)_{k >0}$ such that 
$(\Gamma_{k,n})_{k >0} \equiv  (\Gamma_{k})_{k >0}$ satisfies the assumptions of   Proposition  \ref{Propbeta}. Item (i) then follows by applying Proposition  \ref{Propbeta}. 
\par\ssk
It remains to prove Item (iii). Here we will apply Proposition \ref{PropCriteriaBC} (iii). Thoughout the proof of Item (iii),  $\beta_0 = 1$ by convention. For any positive integer $k$, let 
$S_k = \sum_{j=1}^k  \BBI_{B_j}$ and $E_k = \BBE ( S_k)$.  Since $f$ is convex and $f(0) = 0$, 
$$
f ( (S_k - E_k)/E_n )   \leq (E_k/E_n) f ( (S_k - E_k) / E_k )   
$$
for any  $k$ in $[1,n]$. 
Applying now Inequality (\ref{Laststepbetabis}) in the case $\Gamma_{j,n} = A_j$, we get that 
$$
 \BBE  \bigl(  f ( (S_k - E_k) /E_k  ) \bigr)  \leq E_k^{-1}  \sum_{j=0}^{m-1} 
\int_0^{\beta_j}  ( 1 + 2 Q_{H_k} (s) ) ds  +
2 \int_0^{\beta_m}  Q_{H_k} (s) ds .    
$$
Now, from the  definition of $Q^*$, 
$$
\int_0^u   Q_{H_k} (s) ) ds  \leq u Q^* (u)  \ \text{ for any } u\in ]0,1]\ \text{ and  any }  k>0. 
$$
The three above inequalities ensure that 
\beq  \label{upperboundsupkEfnSk} 
 \sup_{k\leq n } \BBE  \bigl(  f_n ( S_k - E_k  ) \bigr)  \leq E_n^{-1}   \sum_{j\in [0, m-1]}  \beta_j ( 2 Q^* (\beta_j) + 1 ) +
2 \beta_m Q^* ( \beta_m)   . 
\eeq
Let $n_0$ be the smallest integer such that $E_{n_0} \geq 2$.   For $n\geq n_0$, choose $m := m_n = 1 + [E_n]$ in the above inequality.
For this choice of $m_n$, noticing that $Q^* (\beta_j) \geq Q^* (1) = 1$, we get 
\beq
\label{bound1(iii)}
\sum_{n\geq n_0} \frac{\BBP (B_n)}{3 E_n}  \sup_{k\in [1,n]}  \BBE  \bigl(  f_n \bigl( S_k-E_k \bigr) \bigr) \leq 
\sum_{n\geq n_0}  \biggl[ 
 \sum_{ 0\leq j \leq [E_n] }  \beta_j  Q^* ( \beta_j ) \frac{\BBP (B_n)}{E_n^2} 
+ \beta_{m_n} Q^* ( \beta_{m_n} ) \frac{\BBP (B_n)}{E_n}      \biggr]. 
\eeq   
We now bound up the first term on the right-hand side. Clearly 
$$
\sum_{n\geq n_0}  \sum_{ 0\leq j \leq [E_n] }  \beta_j  Q^* ( \beta_j )  E_n^{-2}  \BBP (B_n)   =
\sum_{j \geq 0}  \beta_j  Q^* ( \beta_j ) \sum_{n : E_n \geq j\vee 2 } E_n^{-2}  \BBP (B_n)  . 
$$
Next, noticing that $E_n - E_{n-1} = \BBP (B_n)$, we get that  $\BBP (B_n) / E_n^2 \leq 1/E_{n-1} -1/E_n$.  It follows that 
$$
\sum_{n : E_n \geq j\vee 2 } E_n^{-2}  \BBP (B_n)  \leq 1/E_{n_j-1} \,  ,
$$
where $n_j$ is the smallest integer  such that $E_{n_j} \geq j\vee 2$. Since $E_{n_j - 1} \geq E_{n_j} - 1$, $1/E_{n_j-1} \leq 2 / (j\vee 2)$. 
Hence 
\beq \label{bound2(iii)}
\sum_{n\geq n_0}  \sum_{ 0\leq j \leq [E_n] }  \beta_j  Q^* ( \beta_j )  E_n^{-2}  \BBP (B_n)    \leq 1  +
2 \sum_{j > 0}  j^{-1} \beta_j  Q^* ( \beta_j ) < \infty   
\eeq
under condition (\ref{critstrongBCbeta}). To complete the proof of (iii), it remains to prove that 
\beq \label{laststepth31(iii)}
 \sum_{n\geq n_0}   \beta_{m_n} Q^* ( \beta_{m_n} )  E_n^{-1}  \BBP (B_n)    < \infty 
\eeq 
under condition (\ref{critstrongBCbeta}), where $m_n = 1 + [E_n]$. For any  integer $k \geq 2$, let $I_k$ be the set of integers $n$ such that 
$[E_n] = k$. By definition, $I_k$ is an  interval of $\BBN$. Furthermore, from the fact that $\mu (A_n) \leq 1$, $I_k \not= \emptyset$. Since $\lim_n E_n = \infty$, this interval is finite. Consequently
$$
\sum_{ n \in I_k }  \BBP ( B_n) = E_{\sup I_k}  - E_{\inf I_k - 1}  \leq E_{\sup I_k}  - E_{\inf I_k} +1  \leq 2 . 
$$
Now, recall that $n_0 $ is the first integer such that $E_{n_0} \geq 2$. Consequently $n_0 = \inf I_2$ and 
$$
 \sum_{n\geq n_0}   \beta_{m_n} Q^* ( \beta_{m_n} ) \frac{\BBP (B_n)}{E_n} = \sum_{k\geq 2}
\beta_{k+1} Q^* ( \beta_{k+1} )  \sum_{n\in I_k} \frac{\BBP (B_n)}{E_n}   \leq 2 \sum_{k\geq 2}
k^{-1} \beta_{k+1} Q^* ( \beta_{k+1} )  < \infty \, , 
$$
under condition (\ref{critstrongBCbeta}). This ends the proof of Item (iii).  Theorem \ref{BCbeta} is proved.  $\diamond$

 \subsubsection{Proofs of Theorems \ref{coralphaBC} and \ref{coralphaSBC}  ($\alpha$-mixing case)}

\noindent {\bf Proof of Theorem \ref{coralphaBC}.}  To apply Item (i) of Proposition \ref{PropCriteriaBC}, we shall prove that under \eqref{condregulariryalpha} and \eqref{taschealpha}, there exists a sequence  $(\psi_n)_{n >0}$ of positive integers such that setting $m_n =  \inf \{ k \in \BBN^* \, : \,  \psi_k \geq n\}$, ${\tilde S}_n = \sum_{k=1}^{m_n-1} {\bf 1}_{A_{\psi_k}} (X_{\psi_k}) $ and ${\tilde E}_n = \BBE ({\tilde S}_n ) = \sum_{k=1}^{m_n-1} \mu (A_{\psi_k})$ (so here $g_{j,n}=g_j = \BBI_{A_j} (X_j)$ if $j \in \psi (\mathbb N^*)$ and $0$ otherwise),  we have
\beq \label{CritusualBCbis}
\lim_{N \rightarrow \infty} {\tilde E}_{2^N} = \infty \ \text{ and } \lim_{N \rightarrow \infty}   \bigl(  f \bigl( ( {\tilde S}_{2^N} - {\tilde E}_{2^N})/  {\tilde E}_{2^N}  \bigr) \, \bigr)= 0  \, . 
\eeq
To construct the sequence $\psi= (\psi_n)_{n \geq 1}$, let us make the following considerations. By the second part of  \eqref{taschealpha}, there exists a positive decreasing sequence $(\delta_n)_{n \geq 1} $ such that $\delta_n \rightarrow 0$ , as $n \rightarrow \infty$, and 
\beq \label{taschealphacons1firststep}
\sum_{n \geq 1}  \delta_n \frac{\mu (A_n) }{ \alpha_*^{-1} (  \mu (A_n))} = \infty \, .
\eeq
Now, note that, by the second part of \eqref{condregulariryalpha},   there exist $u_0 >0$ and $\kappa>1$ such that for any $u \in ]0, u_0[$, 
$\alpha_*^{-1} (u/2) \leq \kappa  \alpha_*^{-1} (u)  $.  Hence setting $j_n = \sup \{ j \geq 0 \, : \, \kappa^{-j} \geq \delta_n \}$ and $\varepsilon_n = 2^{-j_n}$, it follows that 
$\alpha_*^{-1} (  \mu (A_n))   \geq \delta_n\alpha_*^{-1} (  \varepsilon_n  \mu (A_n))$,
which combined with \eqref{taschealphacons1firststep} implies that 
\beq \label{taschealphacons1}
\sum_{n \geq 1} \frac{\mu (A_n) }{ \alpha_*^{-1} ( \varepsilon_n \mu (A_n))} = \infty \, .
\eeq

\begin{Definition}  \label{constructionpsin}  
Let $(k_L)_{L \geq 0}$ be the sequence of integers defined by 
\[
k_L=  L \wedge  \lceil \log_2 \alpha_*^{-1} ( \varepsilon_{2^L} \mu (A_{2^L}))  \rceil  \,  , \ \text{ where } \log_2 x = \log (x \vee 1) / \log 2 
\]
and $\lceil x \rceil = \inf \BBZ \cap [x , \infty [$. Set $j_0 =0$ and $j_{L+1} = j_L + 2^{L-k_L}$ for any $L \geq 0$,
Finally, for any $L \geq 0$,  we set $\psi_{j_L} = 2^{L}$ and for any $i=j_L + \ell $ with $\ell \in [1, j_{L+1} - j_L  -1]  \cap \BBN^*$,
$\psi_i = 2^L + \ell 2^{k_L}$. 
\end{Definition} 

Recall the notation, $f_{2^N} (x) = f \bigl( x /  {\tilde E}_{2^N}   \bigr) $. 
Noticing  that ${\tilde S}_{2^N} =  \sum_{k=1}^{j_N-1} {\bf 1}_{A_{\psi_k}} (X_{\psi_k})  $ and recalling that $f(0) =0$, we have
\begin{multline} \label{devissagedirect}
\BBE 
 \bigl(  f \bigl( ( {\tilde S}_{2^N} - {\tilde E}_{2^N})/  {\tilde E}_{2^N}   \bigr) \, \bigr) = \BBE f_{2^N} ({\tilde S}_{2^N}- {\tilde E}_{2^N} )  \\
=  \sum_{L=2}^N  \sum_{\ell =j_{L-1} }^{j_L-1} \Big \{  \BBE f_{2^N} \Big ( \sum_{i=1}^{\ell}  (  {\bf 1}_{A_{\psi_i}}
 (X_{\psi_i}) - \mu(A_{\psi_i}) )  \Big ) -  \BBE f_{2^N} \Big ( \sum_{i=1}^{\ell-1}  (  {\bf 1}_{A_{\psi_i}}
 (X_{\psi_i}) - \mu(A_{\psi_i})  ) \Big )   \Big \}   \, .
\end{multline}
Using Taylor's formula (as to get \eqref{step1beta}) and taking the expectation, we derive 
\begin{multline*}
\BBE f_{2^N} \Big ( \sum_{i=1}^{\ell}  (  {\bf 1}_{A_{\psi_i}}
 (X_{\psi_i}) - \mu(A_{\psi_i})  ) \Big ) -  \BBE f_{2^N} \Big ( \sum_{i=1}^{\ell-1}  (  {\bf 1}_{A_{\psi_i}}
 (X_{\psi_i}) - \mu(A_{\psi_i}) ) \Big )  \\  \leq {\rm Cov} \Big (  f'_{2^N} \Big ( \sum_{i=1}^{\ell-1}  (  {\bf 1}_{A_{\psi_i}}
 (X_{\psi_i}) - \mu(A_{\psi_i}) )  \Big )  ,  {\bf 1}_{A_{\psi_\ell}}
 (X_{\psi_\ell})  \Big ) +  \frac{\mu (A_{\psi_\ell})}{( {\tilde E}_{2^N} )^2}  \, .
\end{multline*}
Since $\Vert f'_{2^N} \Vert_{\infty} \leq 1/  {\tilde E}_{2^N} $, it follows from \eqref{defequivalentalpha} that 
\[
\BBE 
 \bigl(  f \bigl( ( {\tilde S}_{2^N} - {\tilde E}_{2^N})/  {\tilde E}_{2^N}   \bigr) \, \bigr) 
 \leq \sum_{L=2}^N \sum_{\ell=j_{L-1} }^{j_L-1}  \Big \{ \frac{\mu (A_{\psi_\ell})}{( {\tilde E}_{2^N} )^2} +   \frac{ 4 \alpha_{\infty,1} (2^{k_{L-2}})}{{\tilde E}_{2^N} }   \Big \}  \, .
\]
Now, since $j_L - j_{L-1} = 2^{L-k_L}$ and  ${\tilde E}_{2^N}  = \sum_{L=2}^N  \sum_{\ell=j_{L-1} }^{j_L-1}  \mu (A_{\psi_\ell}) $, we get 
\[
\BBE 
 \bigl(  f \bigl( ( {\tilde S}_{2^N} - {\tilde E}_{2^N})/  {\tilde E}_{2^N}   \bigr) \, \bigr) 
 \leq  \frac{1}{ {\tilde E}_{2^N} }+4 \sum_{L=2}^N \alpha_{\infty,1} (2^{k_{L-2}})  \frac{ 2^{L-k_L}}{{\tilde E}_{2^N} }   \, .
\]
Note then that, since $(\mu (A_n))_{n \geq 1}$ is a non-increasing sequence, 
\beq \label{CritusualBCtep0r}
{\tilde E}_{2^N}  = \sum_{L=2}^N  \sum_{\ell=j_{L-1} }^{j_L-1}  \mu (A_{\psi_\ell})   \geq \sum_{L=2}^N  (j_L- j_{L-1} )  \mu (A_{\psi_{j_L}}) = \sum_{L=2}^N  2^{L-k_L} \mu (A_{2^L}) \, .
\eeq
Thus
\[
\BBE 
 \bigl(  f \bigl( ( {\tilde S}_{2^N} - {\tilde E}_{2^N})/  {\tilde E}_{2^N}   \bigr) \, \bigr) 
 \leq  \frac{1}{ {\tilde E}_{2^N} }+  \frac{ 4  \sum_{L=2}^N  \alpha_{\infty,1} ( 2^{k_{L-2}})  2^{L-k_L}}{ \sum_{L=2}^N  2^{L-k_L} \mu (A_{2^L})}   \, .
\]
This shows that \eqref{CritusualBCbis} will be satisfied if 
\beq \label{CritusualBCter}
\lim_{N \rightarrow \infty} {\tilde E}_{2^N} = \infty \ \text{ and } \lim_{L \rightarrow \infty}  
\bigl( \alpha_{\infty,1} ( 2^{k_L}) /  \mu (A_{2^L})\, \bigr) = 0  \, . 
\eeq
Since $(\mu (A_n))_{n \geq 1}$ is a non-increasing sequence, condition \eqref{taschealphacons1} is equivalent to 
\beq \label{taschealphacons2}
\sum_{k \geq 0}  2^k\frac{\mu (A_{2^k}) }{ \alpha_*^{-1} ( \varepsilon_{2^k} \mu (A_{2^k}))} = \infty \, .
\eeq
Together with \eqref{CritusualBCtep0r} and the definition of $2^{k_L}$, \eqref{taschealphacons2} implies the first part of \eqref{CritusualBCter}.  Next, taking into account the definition of $2^{k_L}$, 
\[
  \alpha_{\infty,1} ( 2^{k_L} ) / \mu (A_{2^L})   \leq 
\max \bigl ( C  \varepsilon_{2^L} ,   \alpha_{\infty,1} ( 2^{L}) / \mu (A_{2^L})  \bigr )  \rightarrow 0 \, , \, \text{ as $L \rightarrow \infty $}\, , 
\]
by the first parts of conditions \eqref{condregulariryalpha} and  \eqref{taschealpha}. This ends the proof.  $\diamond$

\medskip

\noindent {\bf Proof of Theorem \ref{coralphaSBC}.} Starting from \eqref{step2beta}, taking into account \eqref{decomp2} and the facts that 
\[
\big | {\rm Cov} \big (  f'_n ( T_{(k-m)_+} ) , {\bf 1}_{A_k} (X_k)   \big )  \big |\leq 4 \alpha_{ \infty,1} (m) /E_n 
\]
and 
\[\big | {\rm Cov} \big (  f'_n ( T_{(k-j)_+} )  - f'_n ( T_{(k-j-1)_+} ), {\bf 1}_{A_k} (X_k)   \big )  \big | 
\leq E_n^{-2} \mu(A_k) \, , 
\]
we infer that, for any positive integer $m$ and any integer $k$ in $[1,n]$,
\[
\BBE  \bigl( f \bigl( (S_k-E_k)/E_n  \bigr) \bigr)  \leq 4 n  \alpha_{ \infty,1} (m) /E_n  +  m/E_n \, .
\]
Item 1. follows by choosing $m=m_n = \eta^{-1} (1/n) $ and  by taking into account Item (ii) of Proposition \ref{PropCriteriaBC}. To prove Item 2., we choose $m =m_n= \alpha_{\infty,1}^{-1} ( u_n E_n/n ) $. Item 2. then follows by taking into account Item (iii) of Proposition \ref{PropCriteriaBC}.  $\diamond$

\subsubsection{Proof of Remark \ref{remreversedtime}} \label{sectionremreversedtime}

To prove that Theorem  \ref{coralphaBC} still holds with $\alpha_{1, \infty} (n)$ replacing  $\alpha_{ \infty, 1} (n)$, it suffices to modify the decomposition \eqref{devissagedirect} as follows: 
\begin{multline*}
\BBE 
 \bigl(  f \bigl( ( {\tilde S}_{2^N} - {\tilde E}_{2^N})/  {\tilde E}_{2^N}   \bigr) \, \bigr) = \BBE f_{2^N} ({\tilde S}_{2^N}- {\tilde E}_{2^N} )  \\
=  \sum_{L=2}^N  \sum_{\ell =j_{L-1} }^{j_L-1} \Big \{  \BBE f_{2^N} \Big ( \sum_{i=\ell}^{j_N-1} (  {\bf 1}_{A_{\psi_i}}
 (X_{\psi_i}) - \mu(A_{\psi_i}) )  \Big ) -  \BBE f_{2^N} \Big ( \sum_{i=\ell +1}^{j_N-1} (  {\bf 1}_{A_{\psi_i}}
 (X_{\psi_i}) - \mu(A_{\psi_i}) ) \Big )   \Big \}   \, .
\end{multline*}
Next, as in the proof of Theorem  \ref{coralphaBC}, we use Taylor's formula and the fact that, by  \eqref{defequivalentalpha},  for any $\ell \in \{j_{L-1}, \ldots, j_{L}-1 \}$, 
\[ {\rm Cov} \Big (  f'_{2^N} \Big ( \sum_{i=\ell +1}^{j_N-1}  (  {\bf 1}_{A_{\psi_i}}
 (X_{\psi_i}) -   \mu(A_{\psi_i}) ) \Big )  ,  {\bf 1}_{A_{\psi_\ell}}
 (X_{\psi_\ell})  \Big )  \leq  \frac{ 4 \alpha_{1, \infty} (2^{k_{L-2}})}{{\tilde E}_{2^N} }    \, .\]
The rest of the proof is unchanged.

To prove that Theorem  \ref{coralphaSBC} still holds with $\alpha_{1, \infty} (n)$ replacing  $\alpha_{ \infty, 1} (n)$, we start by setting
\[
S^*_k =  \sum_{i=k}^n \BBI_{B_{i,n} } \, , \,   E^*_k =  \sum_{i=k}^n  \mu ( \Gamma_{i,n} ) \, , \, T^*_k =S^*_k -E^*_k  \text{ and } \xi_k = T^*_k - T^*_{k+1} \, .
\]
Then, setting $S^*_{n+1}=E^*_{n+1}=0$, instead of \eqref{decomp1}, we write
\[
f_n ( \tilde S_n - \tilde E_n )   =  \sum_{k=1}^n  \bigl(  f_n ( S^*_k - E^*_k) ) - f_n ( S^*_{k+1} - E^*_{k+1} ) \,  \bigr) \, .
\]
By the Taylor integral formula at order $1$, it follows that 
\[
 f_n ( \tilde S_n - \tilde E_n )  =  \sum_{k=1}^n \Big (  f'_n ( T^*_{k+1} )  \xi_k  +   \int_0^1  
\bigl( f'_n ( T^*_{k+1} + t \xi_k ) - f'_n ( T^*_{k+1}) \bigr) \xi_k \, dt   \Big ) \,  . 
\]
Then, instead of \eqref{decomp2}, we use the following decomposition: 
\[
\BBE (  f'_n ( T^*_{k+1} )  \xi_k  )   = \cov ( f'_n ( T^*_{(k +m)_+ } ) , \BBI_{B_{k,n} } )  + 
\sum_{j=1}^{m-1} \cov \bigl( f'_n ( T^*_{(k+j)_+} ) -  f'_n ( T^*_{(k+j+1)_+ } )  ,  \BBI_{B_{k,n} } )  \, .
\]
Hence, the only difference with the proof of Theorem  \ref{coralphaSBC} is the following estimate: 
\[
\big | {\rm Cov} \big (  {\bf 1}_{A_k} (X_k)  ,  f'_n ( T^*_{(k+m)_+} )  \big )  \big |\leq 4 \alpha_{ 1, \infty} (m) /E_n  \, .
\]
This ends the proof of the remark. $\diamond$

\subsection{Proofs of the results of Section \ref{sectioncoefffaible}}

\subsubsection{Proof of Theorem \ref{BCbetafaible}.}

To prove Item  {\rm (i)}, we first apply Proposition \ref{lmacaracterisationBCint}.  Since
$\mu  ( \limsup_n I_n ) >0$,  it follows from that proposition that there exists a sequence $(\Gamma_k)_k$ of intervals
such that $\Gamma_k \subset I_k$,  $\sum_{k>0}  \mu (\Gamma_k)  = \infty$  and 
\begin{equation} \label{infty-norm}
 \sup_{n >0} 
  \left \| \frac{ \sum_{k=1}^n  \BBI_{\Gamma_k} }{\sum_{k=1}^n  \mu(\Gamma_k)}\right \|_{\infty, \mu} < \infty \, ,
\end{equation}
where $\|\cdot\|_{\infty, \mu}$ is the essential supremum norm with respect to $\mu$. 
\par
Let us prove now that $\tilde B_k = \{X_k \in \Gamma_k\}$ is a $ L^1$-Borel-Cantelli sequence. Since $\tilde B_k \subset B_k$,
this will imply that  $( B_k)_{k >0}$ is  a Borel-Cantelli sequence. From (\ref{CritL2BC})  applied to 
$\tilde S_n= \sum_{k=1}^n  \BBI_{\tilde B_k}$,  it is enough to prove that 
\begin{equation}\label{crit-L2}
  \lim_{n \rightarrow \infty} \bigl( {\mathbb E}(\tilde S_n )\bigr)^{-2} \Var (\tilde S_n) =0 \, .
\end{equation}
By stationarity,
\begin{equation}\label{var-eq}
{\mathrm{Var}}(\tilde S_n) =   \sum_{k=1}^n {\mathrm{Var}}(\BBI_{\Gamma_k}(X_0)) +
2 \sum_{k=1}^{n-1} \sum_{j=1}^{n-k} {\mathrm{Cov}} ( \BBI_{\Gamma_k}(X_0), \BBI_{\Gamma_{k+j}}(X_{j}) ) \, .
\end{equation}
Let 
$b_j = \sup_{t \in {\mathbb R}} |{\mathbb E} ( \BBI_{X_j \leq t} | X_0) - {\mathbb P}(X_j \leq t) | $. Clearly, since $\Gamma_\ell$ is  an interval,
\begin{equation}\label{basic-ineq}
|{\mathrm{Cov}} ( \BBI_{\Gamma_k}(X_0), \BBI_{\Gamma_{k+j}}(X_{j}) ) | \leq 2 {\mathbb E}( \BBI_{\Gamma_k}(X_0) b_j) \, .
\end{equation}
Setting $\bar B_{n}= b_1 + \cdots + b_{n-1}$, we infer from \eqref{var-eq}  and \eqref{basic-ineq} that 
\begin{equation} \label{Main-ineq}
{\mathrm{Var}}(\tilde S_n)  \leq {\mathbb E} \Bigl ( (1+ 4\bar B_{n}) \sum_{k=1}^n \BBI_{\Gamma_k}(X_0) \Bigr) \, .
\end{equation}
Since  ${\mathbb E}(\tilde S_n )= \mu(\Gamma_1)+ \cdots + \mu(\Gamma_n)$, we infer from \eqref{Main-ineq} that
\begin{equation}\label{last-step}
\frac{ {\mathrm{Var}}(\tilde S_n)}{\big ( {\mathbb E}(\tilde S_n )\big )^2} \leq  \frac{{\mathbb E} ( 1 +4  \bar B_{n} )}{\sum_{k=1}^n  \mu(\Gamma_k)} 
\left \| \frac{ \sum_{k=1}^n  \BBI_{\Gamma_k} }{\sum_{k=1}^n  \mu(\Gamma_k)}\right \|_{\infty, \mu} 
\leq \frac{ \big (1 + 4 \sum_{k=1}^{n-1} \tilde \beta_{1,1}(k) \big )}{\sum_{k=1}^n  \mu(\Gamma_k)} 
\left \| \frac{ \sum_{k=1}^n  \BBI_{\Gamma_k} }{\sum_{k=1}^n  \mu(\Gamma_k)}\right \|_{\infty, \mu} \, ,
\end{equation}
the last inequality being true because ${\mathbb E}(b_k)= \tilde \beta_{1,1}(k)$. Hence \eqref{crit-L2} follows 
from  \eqref{infty-norm}, \eqref{last-step}, and the fact that $\sum_{k >0}  \tilde \beta_{1,1} (k) < \infty$ and $\sum_{k \geq 1} \mu(\Gamma_k) = + \infty$. The proof of  Item 
{\rm (i)} is complete.

\medskip

We now prove Item {\rm (ii)}. Let $S_n=\sum_{k=1}^n  \BBI_{ B_k}$. Arguing as for {\rm (i)}, it is enough to prove \eqref{crit-L2} with 
$S_n$ instead of $\tilde S_n$. Since the $I_k$ are intervals, the same computations as for {\rm (i)} lead to
\begin{equation} \label{Main-ineq2}
{\mathrm{Var}}( S_j) \leq 
{\mathbb E} \Bigl( (1+4\bar B_{j}) \sum_{k=1}^j \BBI_{I_k}(X_0) \Bigr) 
 \leq {\mathbb E} \Bigl( (1+4\bar B_{n}) \sum_{k=1}^n \BBI_{I_k}(X_0) \Bigr) \, 
\end{equation}
for any $j\leq n$.  Set $\tilde \beta_{1,1} (0)= 1$. Applying H\"older's inequality, we get that, for any $j\leq n$, 
\beq   \label{Holder-ineq}
 {\mathrm{Var}}( S_j)
\leq   \Big\| 1+4 \bar B_n \Big\|_p  \,\,
   \Big \| \sum_{k=1}^n \BBI_{I_k} (X_0) \Big \|_q 
   \leq  4  \Bigl( p \sum_{k=0}^{n-1} (1+k)^{p-1} \tilde \beta_{1,1}(k) \Bigr)^{1/p}
    \Big \| \sum_{k=1}^n \BBI_{I_k} (X_0) \Big \|_q \, ,
\eeq
(the last inequality follows from Remark 1.6 and Inequality (C.5) in \cite{Rio17}). 
Consequently
$$
\frac{ {\mathrm{Var}}( S_n)}{\big ( {\mathbb E}( S_n )\big )^2}  \leq \frac {4   
\bigl( p  \sum_{k=0}^{n -1} (1+k)^{p-1} \tilde \beta_{1,1}(k) \bigr)^{1/p}}{E_n}
    \left \| \frac{\sum_{k=1}^n \BBI_{I_k} (X_0)}{E_n} \right \|_q \, .
$$
Hence  Item {\rm (ii)} follows via  (\ref{CritL2BC}). In addition, 
Item {\rm (iii)} follows  from \eqref{Holder-ineq} by applying (\ref{CritstrongBCvariance}). 

\medskip 
To prove {\rm (iv)} and {\rm (v)}, we start from \eqref{Main-ineq2}, and we get that, for any $j\leq n$, 
\begin{equation}\label{Holder-ineq-bis}
 \frac{ {\mathrm{Var}}( S_j)}{\big ( {\mathbb E}( S_n )\big )^2}
\leq  \frac{\|1+4 \bar B_n \|_\infty}{E_n^2}  
    \sum_{k=1}^n \mu(I_k)  
   \leq  \frac{1}{E_n} \Bigl(1 + 4\sum_{k=1}^{n-1}  \tilde \varphi_{1,1}(k) \Bigr)  \, .
\end{equation}
Then  {\rm (iv)} follows from \eqref{Holder-ineq-bis} with $j=n$ and  (\ref{CritL2BC}) and 
 {\rm (v)}  from \eqref{Holder-ineq-bis} and (\ref{CritstrongBCvariance}).  $\diamond$

\subsubsection{Proof of Lemma \ref{easy}}

We consider  the natural coupling
$$
X_k^*= \sum_{i=0}^{k-1} 2^{-i} \varepsilon_{k-i}+\sum_{i \geq k} 2^{-i}\varepsilon'_{k-i} \, ,
$$
where $(\varepsilon'_i)_{i \in {\mathbb Z}}$ is an
independent copy of $(\varepsilon_i)_{i \in {\mathbb Z}}$. Note that $X_k^*$ distributed as $X_k$ and independent of 
$X_0$. 
Let then 
$$
  \delta(k)=  {\mathbb E}\left(\min(|X_k-X_k^*|, 1)\right) \, .
$$
We  first give a bound on $\delta(k)$. By definition 
$$
\delta(k) \leq {\mathbb E}
\Bigl(\min \Bigl( \, \Big|\sum_{i \geq k} 2^{-i} (\varepsilon_{k-i}-\varepsilon'_{k-i} )\Big|, 
1 \Bigl)\, \Bigr) \, .
$$
By sub-additivity and stationarity,
$$
\delta(k) \leq 
\sum_{i \geq k} {\mathbb E}\left( \min(2^{-i} |\varepsilon_0-\varepsilon_0'|, 1)\right ) \, .
$$
Hence
\begin{equation*}
\delta(k) \leq 
 \sum_{i \geq k} 2^{-i} {\mathbb E}\Bigl(  |\varepsilon_0-\varepsilon_0'| {\bf{1}}_{|\varepsilon_0-\varepsilon_0'|\leq 2^{i/2}} \Bigr) 
+\sum_{i \geq k}
{\mathbb P}\left ( |\varepsilon_0-\varepsilon_0'| > 2^{i/2}\right )
\end{equation*}
and, consequently,
\begin{equation}\label{1delta}
\delta(k) \leq  \kappa 2^{-k/2} + {\mathbb E}\bigl(  \bigl(   2 (\log 2)^{-1}  \log |\varepsilon_0-\varepsilon_0'|
-k\bigr)_+\bigr) 
\end{equation}
with $\kappa= 1/(1-2^{-1/2})$. This gives the upper bound 
$$
\delta(k) \leq  K 2^{-k/2} + K{\mathbb E}\bigl(
(  \log |\varepsilon_0-\varepsilon_0'|
{ \bf 1}_{\log |\varepsilon_0-\varepsilon_0'|>k \log(\sqrt 2)}\bigr) \, .
$$
Now, if \eqref{weakp} holds,
 $$
\sup_{t >1} t^{p-1} {\mathbb E}\bigl(
 \log |\varepsilon_0-\varepsilon_0'|
{ \bf 1}_{\log |\varepsilon_0-\varepsilon_0'|>t }\bigr) < 
\infty \, ,
$$ 
and it follows then easily from \eqref{1delta} that there exists some positive constant $B$ such that
\begin{equation}\label{2delta}
  \delta (k) \leq Bk^{1-p}  \ \text{ for any }  k\geq 1.  
\end{equation}
\par
Now let $F_\mu$ be the distribution function of $\mu$.  By Lemma 2, Item 2.  in \cite{DP}, for any $y \in [0,1]$
\beq \label{majorbetatilde}
\tilde \beta_{1,1}(k) \leq  y +{\mathbb P}(|F_\mu (X_k)-F_\mu (X_k^*)|>y) 
\eeq
Since $\mu$ has a bounded density, $F_\mu$ is Lipshitz.  Moreover $|F_ \mu (X_k)-F_\mu (X_k^*)|\leq 1$. Hence  
\beq\label{diffFmu}
|F_\mu (X_k)-F_\mu (X_k^*)| \leq  A \min(1 ,  |X_k-X_k^*| ) \ \text{ for some  constant } A \geq 1.  
\eeq
Now, by (\ref{majorbetatilde}), (\ref{diffFmu}) and the  Markov inequality, 
$\tilde \beta_{1,1}(k) \leq  y + A \delta(k) / y$ for any positive $y$. Consequently 
$\tilde \beta_{1,1}(k) \leq 2 \sqrt{A \delta (k)}$.  The conclusion of Lemma \ref{easy} follows then from  \eqref{2delta}. $\diamond$

\subsubsection{Proof of Lemma \ref{notsoeasy}}

We first note that, for any function $g$ in $L^2(\lambda)$, one has 
\begin{equation}\label{depart}
  K^n(g)(x)- \lambda (g) = \sum_{k \in {\mathbb Z}^*}  (\cos (2 \pi k a ))^n \hat g(k) \exp(2 i \pi k x) \, ,
\end{equation}
where $(\hat g(k))_{k \in {\mathbb Z}}$ are the Fourier coefficients of $g$. 

Next, we need to approximate the function ${\bf 1}_{[0,t]}$ by smooth functions. To do this, we start  from  an infinitely differentiable density 
$ \ell$ supported in $[0,1]$, and we define 
$$g_1(x)=  \left (\int_0^x \ell(t) dt  \right ) {\bf 1}_{[0,1]}(x)  \quad \text{and} \quad g_2(x)= (1-g_1(x)) {\bf 1}_{[0,1]}(x) \, .$$
Now, for $0<h <1/4$,  $t \in [2h, 1-h]$ and $ x \in [0,1]$, we have
$$
     f_{t,h}^-(x) \leq {\bf 1}_{[0,t]}(x) \leq f_{t,h}^+(x)\, ,
$$
where
\begin{align*}
  f_{t,h}^+(x) &= {\bf 1}_{[0,t]}(x) +g_2( (x-t)/h  ) + g_1 ( (x+h-1)/h  ) \\
  f_{t,h}^-(x) &= {\bf 1}_{[h,t-h]}(x) +g_2 ( (x+h-t)/h ) + g_1 (x/h)  \, .
\end{align*}
Hence, for $t \in [2h, 1-h]$
\begin{equation}\label{B1}
K^n(f_{t,h}^-) -\lambda (f_{t,h}^-) - 2h \leq K^n( {\bf 1}_{[0,t]}) -t \leq K^n(f_{t,h}^+) -\lambda (f_{t,h}^+) + 2h \, .
\end{equation}
On the other hand
\begin{equation}\label{B2}
\Big \|\sup_{t \in [0,2h]} | K^n({\bf 1}_{[0,t]})-t | \Big  \|_{1}  \leq 4h  \quad \text{and}  \quad 
\Big \|\sup_{t \in [1-h,1]} | K^n({\bf 1}_{[0,t]})-t | \Big  \|_{1}  \leq 2h \, .
\end{equation}
From \eqref{B1} and \eqref{B2}, we get
\begin{equation}\label{mainineq}
\Big \|\sup_{t \in [0,1]} | K^n({\bf 1}_{[0,t]})-t | \Big  \|_{1}  \leq 10h + 
\Big \|\sup_{f\in {\mathcal F}_h} | K^n(f)-\lambda(f) | \Big  \|_{1} 
\end{equation}
where ${\mathcal F}_h= \{ f_{t,h}^+, f_{t,h}^-, t \in [2h, 1-h] \}$. 

Note that the functions belonging to ${\mathcal F}_h$ are infinitely differentiable, so that  one can easily find some upper bounds 
on their Fourier coefficients. More precisely, by two elementary integrations by parts, we obtain that there exist a positive constant $C$ such that, 
for any $f \in {\mathcal F}_h$, 
\begin{equation} \label{F1}
| \hat f(k) | \leq \frac {C}{h (|k|+1)^2}\, .
\end{equation}

From \eqref{depart} and \eqref{F1}, we get that
\begin{equation}\label{phibound}
\sup_{f\in {\mathcal F}_h} \Vert K^n (f) - \lambda (f) \Vert_{\infty, \lambda } \leq 
\frac C h  \sum_{k \in {\mathbb Z}^*} |k|^{-2} |\cos (2 \pi k a )|^n \, .
\end{equation}
Take $\beta \in (0, 1/2)$. By the properties of the Gamma function there exists a positive constant $K$ such that, 
$$
\sum_{n \geq 1}  \frac{1}{n^{\beta}} |\cos(2\pi k a)|^n  \leq   \frac{K}{(1 - | \cos (2 \pi ka) | )^{1- \beta}} \, .
$$
Since $(1-|\cos(\pi u)|) \geq \pi
(d(u, {\mathbf{Z}}))^2$, we derive that
$$
\sum_{n \geq 1} \frac{1}{n^{\beta}} \sum_{k \in {\mathbf{Z}}^*} |k|^{-2}
|\cos(2\pi k a)|^n  \leq   \frac{K}{\pi^{1- \beta}} \sum_{k \in {\mathbf{Z}}^*} \frac{|k|^{-2}} {(d(2ka, {\mathbf{Z}}) )^{2-2 \beta}} \, .
$$
Note that, if $a$ is badly approximable by
rationals in the weak sense, then so is $2a$.
Therefore if $a$ satisfies \eqref{badly}, proceeding as in the proof of Lemma 5.1 in \cite{DR},  we get that, for any $\eta>0$,
\begin{equation*}
\sum_{k=2^N}^{2^{N+1}-1} \frac{1}{(d(2ka, {\mathbf{Z}}) )^{2- 2\beta}} =
{\cal O} (2^{(2-2\beta)N (1+ \eta)})  \, .
\end{equation*}
Therefore, since $\beta \in (0,1/2)$, taking $\eta$ close enough to 0, we get 
\begin{equation}\label{nearend}
\sum_{n \geq 1} \frac{1}{n^{\beta}}\sum_{k \in {\mathbf{Z}}^*} |k|^{-2} 
|\cos(2\pi k a)|^n   
 \ll
\sum_{N \geq 0} 2^{(2-2\beta)N (1+ \eta)}  \max_{2^N \leq k \leq
2^{N+1}} |k|^{-2}< \infty \, .
\end{equation}
From \eqref{phibound} and \eqref{nearend}, for any $c$ in $(0,1)$ there exists a constant $B$ such that
\begin{equation} \label{nearend2}
\sup_{f\in {\mathcal F}_h} \Vert K^n (f) - \lambda (f) \Vert_{\infty, \lambda }  \leq B n^{-c} h^{-1} \, .
\end{equation}
From \eqref{mainineq} and \eqref{nearend2}, we infer that, for any $c$ in $(0,1)$ there exists a  constant $\kappa$ such that
$$
\Big \|\sup_{t \in [0,1]} | K^n({\bf 1}_{[0,t]})-t | \Big  \|_{1}  \leq \kappa \bigl( h + n^{-c} h^{-1}   \bigr) \, .
$$
Taking $h = n^{c/2}$ in the above inequality, we then get Lemma \ref{notsoeasy}.  $\diamond$

\subsubsection{Proof of Corollary \ref{RWwithdrift}}

The first part of Corollary \ref{RWwithdrift} follows immediately from Lemma \ref{notsoeasy}  and Theorem \ref{thmalphafaible} applied to 
$\XZ$ and the sequence $(J_n)$ of intervals  on the circle defined by $J_n = [nt, nt + n^{-\delta} ]$.  In order to prove the second part, we will apply Theorem \ref{BCbetafaible}(iii) to the sequence $\XZ$.  The main step is to prove that 
\beq
\label{CRW1} 
\sup_{n >0}  \frac{1} {E_n}   \Bigl \| \sum_{k=1}^n \BBI_{J_k} (X_0) \Bigr \|_\infty < \infty . 
\eeq
Now  $E_n \sim n^{1-\delta} / (1-\delta) $ as $n\rightarrow \infty$. Therefrom one can easily see that 
 (\ref{CRW1}) follows from the inequality below: for some positive constant $c_0$, 
 \beq
\label{CRW2} 
\sum_{k=m+1}^{2m} \BBI_{J_k}  \leq c_0 m^{1-\delta} \text{ for any integer } m>0. 
\eeq
Now 
$\sum_{k=m+1}^{2m} \BBI_{J_k}  (x) \leq \sum_{k=m+1}^{2m} \BBI_{x - m^{-\delta} \leq kt \leq x}$.
Furthermore, if $t$ is badly approximable, 
then, from (\ref{badly}) with $\eps = 0$,  $d ( kt , lt)  = d ( t (l-k) , \BBZ ) \geq c (l-k)^{-1}  \geq c/m$
for any  $(k,l)$ such that $m< k <l \leq 2m$, which ensures that  
$\sum_{k=m+1}^{2m} \BBI_{x - m^{-\delta} \leq kt \leq x} \leq 1 + c^{-1} m^{1-\delta}$ for any $x$. 
This inequality and the above facts imply (\ref{CRW2}) and, consequently, (\ref{CRW1}).  Now Corollary \ref{RWwithdrift} follows easily from
Lemma \ref{notsoeasy},  (\ref{CRW1}) and Theorem \ref{BCbetafaible}(iii)  $\diamond$

\subsection{Proofs of the results of Section \ref{HarrisMC}}

\subsubsection{Proof of Theorem \ref{BCHarris}.} Recall that, for any Polish space $E$, there exists a one to one bimeasurable mapping from $E$ onto a Borel subset of $[0,1]$.
Consequently we may assume without loss of generality that $E= [0,1]$.  
We define the Markov chain and the renewal process  in the same way as in Subsection 9.3 in Rio (2017). Let
 $(U_i, \eps_i)_{i\geq 0}$ be a sequence of independent random variables with the uniform law over $[0,1]^2$ and  $\zeta_0$ be a random variable with law $\mu$ 
 independent of $(U_i,\eps_i)_{i\geq 0}$. 
Let $(\xi_k)_{k>0}$ be a sequence of independent random variables with law $\nu$. Suppose furthermore that this sequence $(\xi_k)_{k>0}$
is independent of the $\sigma$-field generated by $\zeta_0$ and $(U_i,\eps_i)_{i\geq 0}$. Define the stochastic kernel $Q_1$ by 
\beq \label{defQ1}
Q_1 (x, A) = (1- s(x) )^{-1} ( P (x.A) - s(x) \nu (A) ) \ \text{if }s(x) < 1 \ \text{and } Q_1 (x,A) = \nu (A) \ \hbox{if } s(x) = 1  
\eeq
and the conditional distribution function $G_x$ by 
\beq \label{defcdfG}
G_x  (t) = Q_1 ( x,  \,  ] - \infty , t] \,  )  \ \text{ for any } (x,t) \in [0,1] \times [0,1] . 
\eeq
Define the sequence $(X_n)_{n>0}$ by induction in the following way:  $X_0 = \zeta_0$ and 
\beq\label{defXn}
X_{n+1}  = \xi_{n+1} \ \text{if } s(X_n) \geq U_n \ \text{and } X_{n+1} = G_{X_n}^{-1}  (\eps_n)\ \hbox{if } s(X_n) < U_n .  
\eeq
Then the sequence $(X_n)_{n\geq 0}$ is a  Markov chain with kernel $P$ and initial law $\mu$. The incidence process $(\eta_n)_{n\geq 0}$
is defined by $\eta_n = \BBI_{U_n \leq s(X_n)}$ and the renewal times $(T_k)_{k\geq 0}$ by 
\beq \label{defTk}
T_k = 1 + \inf \{  j\geq 0 : \eta_0 + \cdots + \eta_j = k+1 \}   . 
\eeq
We also set $\tau_j = T_{j+1} - T_j$ for any $j\geq 0$. Under the assumptions of Theorem \ref{BCHarris},  $(\tau_j)_{j\geq 0}$ is 
a sequence of integrable, independent and indentically distributed random variables. Note also that \eqref{Pitmanfinitemass} 
implies that $T_0 < \infty$ almost surely (see Rio (2017), Subsection 9.3). Hence, by the strong law of large numbers, 
\beq\label{SLLNTk}
\lim_{k \rightarrow \infty}  (T_k/k)  = \BBE ( \tau_1 )   \ \hbox{a.s.}
\eeq
Let  $m$ be a positive integer such that $m> \BBE (\tau_1)$.  Then there exists some random integer $k_0$ such that
$T_k \leq km$  for any $k\geq k_0$.
Since the sequence of sets $(A_j)_{j>0}$ is non-increasing, it follows that 
$\BBI_{A_{T_k} }  \geq \BBI_{A_{km} }$ for any $k\geq k_0$.
Furthermore  
\beq
 \sum_{k>0}  \BBI_{A_k} (X_k)  \geq \sum_{k>0} \BBI_{A_{T_k} } ( X_{T_k} ) .
\eeq
Consequently,  if $\sum_{k>0} \BBI_{A_{km} }  ( X_{T_k} )  = \infty$ a.s.,  then  a.s.  $\sum_{k>0}  \BBI_{A_k} (X_k)  = \infty$. Now, from the construction of the Markov chain, the random variables $(X_{T_k} )_{k>0}$ are iid with law $\nu$.
Next, since  the sequence of sets $(A_j)_{j>0}$ is non-increasing and $\sum_k \nu (A_k) = \infty$, the series 
$\sum_k \nu ( A_{km} )$ is divergent. Hence, by the second Borel-Cantelli lemma for sequences of independent events, 
$\sum_{k>0} \BBI_{A_{km} }  ( X_{T_k} )  = \infty$ a.s., which completes the proof of Theorem \ref{BCHarris}.  $\diamond$

\subsubsection{Proof of Theorem \ref{BCHarrisconverse}.}   From Lemma 9.3 in Rio (2017), the stochastic kernel $P$ is 
irreducible, aperiodic and positively recurrent. 
Furthermore 
$$
\mu =   \biggl( \int_E  {1\over s(x) }  d\nu (x)  \biggr)^{-1}  \frac{1}{s(x)}  \, \nu 
$$
is the unique invariant law under $P$. Now, let $(X_i)_{i \in \BBN}$ denote the strictly stationary Markov chain with kernel $P$. 
Define the renewal times $T_k$ as in \eqref{defTk}. Then the random variables $(X_{T_k})_{k>0}$ are iid with law $\nu$. 
Since $\sum_{k>0} \nu (A_k) < \infty$, it follows that $\sum_{k>0}  \BBI_{A_k}  ( X_{T_k} ) < \infty$ almost surely. 
Now $T_k \geq k$, from which $A_{T_k} \subset A_k$. Hence 
\beq\label{BCTk}
\BBP  ( X_{T_k} \in A_{T_k} \ \text{infinitely often }  ) = 0 . 
\eeq
Since $Q_1 (x, .) = \delta_x$, $X_m = X_{T_k}$ for any $m$ in $[T_k , T_{k+1} [$. Furthermore $A_m\subset A_{T_k}$ for 
any $m\geq T_k$. Consequently, if $X_{T_k}$ does not belong
to $A_{T_k}$, then, for any $m$ in $[T_k , T_{k+1} [$, $X_m$ does not belong $A_m$. Now (\ref{BCTk}) and the above fact imply
Theorem \ref{BCHarrisconverse}.  $\diamond$

\appendix

\section{Uniform integrability}
\label{uniformint}

\setcounter{equation}{0}

In this section, we recall the definition of the uniform integrability and we give a criterion for the uniform integrability of a family
$(Z_i)_{i\in I}$ of nonnegative random variables. We first recall the usual definition of uniform integrability, as given in Billingsley (1999).

\begin{Definition}  A family $(Z_i)_{i\in I}$ of nonnegative random variable is said to be uniformly integrable if
$\lim_{M\rightarrow +\infty}  \sup_{i\in I}  \BBE \bigl( Z_i \BBI_{Z_i > M}  \bigr)   = 0$. 
\end{Definition}

Below we give a proposition, which provides a more convenient criterion. In order to state this proposition, we need to introduce some
quantile function.

\begin{Notation}   Let $Z$ be a real-valued random variable and $H_Z$ be the tail function of $Z$, defined by $H_Z (t) = \BBP (Z>t)$ 
for any real $t$. We denote by $Q_Z$ the cadlag inverse of $H_Z$. 
\end{Notation}

\begin{Proposition}   \label{propquantileuniformint}  A family $(Z_i)_{i\in I}$ of nonnegative random variables is  uniformly integrable
if and only if 
\beq\label{critquantileuniformint} 
\lim_{\eps \searrow 0}  \sup_{i\in I}  \int_0^\eps Q_{Z_i} (u) du = 0 .
\eeq
\end{Proposition}

\par\no
{\bf Proof.}  Assume that the family $(Z_i)_{i\in I}$ is uniformly integrable. 
Let $U$ be a random variable with uniform distribution over $[0,1]$.
Since $Q_{Z_i} (U)$ has the same distribution as $Z_i$,  
$$
\sup_{i\in I} \int_0^\eps Q_{Z_i} (u) du \leq   M\eps  +   \sup_{i\in I}  \BBE \bigl( Z_i \BBI_{Z_i > M}  \bigr) . 
 $$
Choosing $M = \eps^{-1/2}$ in the above inequality, we then get (\ref{critquantileuniformint}). 
Conversely, assume that condition  (\ref{critquantileuniformint}) holds true. Then one can easily prove that 
$A: = \sup_{i\in I} \BBE (Z_i) < \infty$.  It follows that $\BBP ( Z_i > A/\eps )  \leq \eps $, which ensures that 
$Q_{Z_i} (\eps) \leq A/\eps$. Consequently, for any $i\in I$,  
$$
\BBE \bigl( Z_i \BBI_{Z_i > A/\eps}  \bigr)   = \int_0^{\eps} Q_{Z_i} (u) \BBI_{Q_{Z_i} (u) > A/\eps} \, du \leq \int_0^{\eps} Q_{Z_i} (u) du, 
$$
which implies the uniform integrability of $(Z_i)_{i\in I}$.  $\diamond$

\section{Criteria under pairwise correlation conditions }
\label{pairwisecoeff}

\begin{Proposition}   \label{propPC}   Let $(B_k)_{k>0}$ be a sequence of events in 
$\EP$ such that $\BBP (B_1) >0$ and $\sum_{k>0} \BBP (B_k) = \infty$.  Set $E_n = \sum_{k=1}^n  \BBP ( B_k) $. 
Assume that there exist a non-increasing sequence $(\gamma_n)_n $ of reals in $[0,1]$ 
and sequences $(\alpha_n)_n $ and $(\varphi_n)_n $ of reals in $[0,1]$ such that for any integers $k$ and $n$, 
\[
\big | \BBP ( B_k \cap B_{k+n} ) -   \BBP ( B_k)   \BBP ( B_{k+n} )  \big | \leq  \gamma_n \BBP ( B_k) \BBP ( B_{k+n} ) +  \varphi_n \big (  \BBP ( B_k) +  \BBP ( B_{k+n} ) \big ) + \alpha_n \, .
\]
\par\ssk\no
{\rm (i)}  Assume that 
\begin{equation} \label{eqpropPC1}  
\gamma_n \rightarrow 0 \  , \  E_n^{-1}   \sum_{k=1}^n \varphi_k \rightarrow 0   \,  \text{ and } \, E_n^{-2}   \sum_{k=1}^n \sum_{j=1}^k \min (\alpha_j ,\BBP (B_k) )  \rightarrow 0 \, ,  \text{ as } n \rightarrow \infty \, .
\end{equation}
Then $(B_k)_{k>0}$ is a $L^1$ Borel-Cantelli sequence. 
\par\ssk\no
{\rm (ii)}   Assume that 
\begin{equation} \label{eqpropPC2}
  \sum_{k\geq 1}  \frac{ \gamma_k }{ k } < \infty  \  , \ \sum_{k\geq 1}  \frac{ \varphi_k}{E_k } < \infty  \,  \text{ and } \, \sum_{k\geq 1}  E_k^{-2}    \sum_{ j \in [1, k] } \min (\alpha_j ,\BBP (B_k) )  < \infty \, .
\end{equation}
Then $(B_k)_{k>0}$ is a strongly Borel-Cantelli sequence. 
\end{Proposition}
\begin{Remark} \label{RemarkpropBC}
If $\alpha_n = {\cal O} (n^{-a}) $ with $ a \in ]0,1[$, then $ \sum_{j=1}^k \alpha_j = {\cal O} (k^{1-a})$. Hence the third condition in \eqref{eqpropPC1} holds as soon as $n^{ -1 + (a/2)} E_n \rightarrow \infty$. On the other hand, the third condition in \eqref{eqpropPC2} holds as soon as 
$ \sum_{n \geq 1} n^{1-a} E_n^{-2} < \infty$ (note that this latter condition is satisfied when $ n^{-1+a/2} (\log n)^{- ( 1/2 + \varepsilon) }E_n \rightarrow \infty$ for some $\varepsilon >0$). 
\par
If $\alpha_n = {\cal O} (n^{-1}) $  then $\sum_{j=1}^k \alpha_j = {\cal O} (\log k)$. Hence the third condition in \eqref{eqpropPC1} holds as soon as $E_n (n \log n )^{-1/2} \rightarrow \infty$. On the other hand, the third condition in \eqref{eqpropPC2} holds as soon as $ \sum_{n \geq 1}  (n \log n ) /E_n^2 < \infty$ (note that this latter condition is satisfied when $ n^{-1/2}(\log n)^{- ( 1 + \varepsilon )} E_n \rightarrow \infty$ for some $\varepsilon >0$).
\par
If $\alpha_n = {\cal O} (n^{-a}) $ with $ a>1$, then  
$ \sum_{j = 1}^\infty  \min (\alpha_j ,\BBP (B_k) )  =   {\cal O} \bigl( \BBP (B_k)^{1-1/a} \bigr)$.  
Hence the third condition in \eqref{eqpropPC1} holds if $n^{-1/(a+1)} E_n \rightarrow \infty$ 
(use the fact that $\sum_{k=1}^n \BBP (B_k)^{1- 1/a} \leq n ( E_n/n)^{{1-1/a}}$). Next, the third condition in \eqref{eqpropPC2} holds as soon as $ \sum_{n \geq 1} E_n^{-2} \BBP (B_n)^{1-1/a}   < \infty$ (note that this latter condition is satisfied when $ \BBP (B_n) 
\geq n^{-a/(a+1) } ( \log n)^{a/(a+1) + \varepsilon }$  for some $\varepsilon >0$). 
\par
If $\alpha_n = {\cal O} ( a^n )$ with $ a \in ]0,1[$ then 
$  \sum_{j= 1}^\infty  \min (\alpha_j ,\BBP (B_k) )     = O \bigl( \BBP (B_k) \log \big (  e /  \BBP (B_k) \bigr) $. 
Hence  the third condition in \eqref{eqpropPC1} holds as soon as $n \,  \BBP (B_n)  \rightarrow \infty$. 
On the other hand, the third condition in \eqref{eqpropPC2} holds as soon as  $ \BBP (B_n) 
\geq n^{-1 } ( \log n)^{ \varepsilon }$ for some $\varepsilon >0$. 
\end{Remark}

\par\no
{\bf Proof of Proposition \ref{propPC}.}  
Note that 
\begin{multline}  \label{p1propPC17-dec}
\max_{ k \leq n}{\rm Var} S_k \leq E_n + 2 \sum_{i=1}^{n-1}\sum_{j=i+1}^n   \Big ( \gamma_ {j-i} \BBP ( B_i) \BBP ( B_{j} ) +  \varphi_ {j-i} \big (  \BBP ( B_i) +  \BBP ( B_{j} ) \big ) +  (    \BBP ( B_j) \wedge \alpha_ {j-i}  ) \Big ) \\
\leq E_n  \big ( 1+ 4 \sum_{k=1}^{n-1}    \varphi_ {k}  \big )  + 2  \sum_{j=2}^n  \sum_{k=1}^{j-1} 
 (    \BBP ( B_j)  \wedge\alpha_ k )  +  2 \sum_{i=1}^{n-1}\sum_{j=i+1}^n  \gamma_ {j-i} \BBP ( B_i) \BBP ( B_{j} )   \, .
\end{multline}
Moreover, for any positive integer $m$,
\begin{align}  \label{p2propPC17-dec}
\sum_{i=1}^{n-1}\sum_{j=i+1}^n  \gamma_ {j-i} \BBP ( B_i) \BBP ( B_{j} )  & \leq  \sum_{i=1}^{n-1}\sum_{j=i+1}^{(i+m-1) \wedge n} \BBP ( B_i) +  \gamma_m \sum_{i=1}^{n-1}\sum_{j=i+m}^n   \BBP ( B_i) \BBP ( B_{j} )  \nonumber  \\
& \leq m E_{n} + \gamma_m  E_{n-1} E_n \, .
\end{align}
Now, from \eqref{p1propPC17-dec} and \eqref{p2propPC17-dec}, one easily infers that criteria   \eqref{CritL2BC} holds true under 
 \eqref{eqpropPC1}, which proves Item (i) of Proposition \ref{propPC}.    

To prove Item (ii), we shall apply criteria \eqref{CritstrongBCvariance}. Starting from \eqref{p1propPC17-dec} and using the facts that $\sum_{j \geq k} E_j^{-2} \BBP (B_j) \leq E_{k-1}^{-1}$ and $\sum_{n\geq j} E_n^{-3} \BBP ( B_n) \leq E_{j-1}^{-2} $, we get that 
\[
\sum_{n \geq 2}  \frac{\BBP (B_n)}{4E_n^3}  \max_{ k \leq n}{\rm Var} S_k  \leq \frac{2}{E_1} +   \sum_{k\geq 2}  \frac{ \varphi_k}{E_{k-1} }  +    \sum_{j \geq 2}  \sum_{k=1}^{j-1}  \frac{\BBP ( B_j) \wedge \alpha_ {k}  } {E_{j-1}^2}   
+  \sum_{n \geq 2}  \frac{\BBP (B_n)}{E_n^2}  \min_{ m \geq 1} \big ( m  + \gamma_m E_{n-1} \big ) .
\]
By the second and the third conditions in \eqref{eqpropPC2}, it follows that \eqref{CritstrongBCvariance} will be satisfied if
\begin{equation} \label{conditionpourgamma}
 \sum_{n \geq 2}  \frac{\BBP (B_n)}{E_n^2}  \min_{ m \geq 1} \big ( m  + \gamma_m E_{n-1} \big )  < \infty \, .
\end{equation}
Define the function $\psi : [1, \infty [ \mapsto [0, \infty [$ by  $\psi (x) = ( \gamma_{[x]} / [x] )$ and let 
$\psi^{-1}$ denote the cadlag generalized inverse function of $\psi$.  Let  
$m_n = \psi^{-1} (E_{n-1}^{-1})$. Then $m_n \geq 1$ and 
\[ 
\min_{ m \geq 1} \big ( m  + \gamma_m E_{n-1} \big )  \leq m_n + m_n \psi ( m_n) E_{n-1} \leq 2 m_n ,
\]
since $\psi ( \psi^{-1} (x) ) \leq x$.  Using the fact that $ x \mapsto \psi^{-1} (1/x) $ is non-decreasing, it follows that 
\begin{multline*}
\sum_{n \geq 2}  \frac{\BBP (B_n)}{E_n^2}  \min_{ m \geq 1} \big ( m  + \gamma_m E_{n-1} \big )    \leq  2  \sum_{n \geq 2}  \frac{\BBP (B_n)}{E_n^2}  \psi^{-1} (E_{n-1}^{-1})  \\
 \leq 2 \int_{0}^{1/E_1}  \psi^{-1} ( x) dx =  2  \int_{ \psi^{-1} (1/E_1)}^{\infty}  \psi ( y) dy +  \frac{2}{E_1}  \psi^{-1} (1/E_1) \, , 
\end{multline*}
which is finite under the first part of condition \eqref{eqpropPC2}. This ends the proof of Item (ii).  $\diamond$


\end{document}